\documentclass[11pt]{article}

\usepackage{latexsym}
\usepackage{amsfonts}
\usepackage{amssymb}
\usepackage{amsmath}
\usepackage{mathrsfs}
\usepackage[all]{xy}
\usepackage{epic}
\usepackage{eepic}

\newcommand{\bldr}{\mathbf{r}}
\newcommand{\blds}{\mathbf{s}}

\newcommand{\bldt}{\mathbf{t}}

\newcommand{\bldz}{\mathbf{z}}

\newcommand{\scrc}{\mathscr{C}}
\newcommand{\scrf}{\mathscr{F}}
\newcommand{\be}{\begin{equation}}
\newcommand{\ee}{\end{equation}}
\newcommand{\bea}{\begin{eqnarray}}
\newcommand{\eea}{\end{eqnarray}}
\newcommand{\bean}{\begin{eqnarray*}}
\newcommand{\eean}{\end{eqnarray*}}
\newcommand{\brray}{\begin{array}}
\newcommand{\erray}{\end{array}}

\newtheorem{dfn}{Definition}[section]
\newtheorem{thm}[dfn]{Theorem}
\newtheorem{lmma}[dfn]{Lemma}
\newtheorem{ppsn}[dfn]{Proposition}
\newtheorem{crlre}[dfn]{Corollary}
\newtheorem{xmpl}[dfn]{Example}
\newtheorem{rmrk}[dfn]{Remark}

\newcommand{\bdfn}{\begin{dfn}}
\newcommand{\bthm}{\begin{thm}}
\newcommand{\blmma}{\begin{lmma}}
\newcommand{\bppsn}{\begin{ppsn}}
\newcommand{\bcrlre}{\begin{crlre}}
\newcommand{\bxmpl}{\begin{xmpl}}
\newcommand{\brmrk}{\begin{rmrk}}

\newcommand{\edfn}{\end{dfn}}
\newcommand{\ethm}{\end{thm}}
\newcommand{\elmma}{\end{lmma}}
\newcommand{\eppsn}{\end{ppsn}}
\newcommand{\ecrlre}{\end{crlre}}
\newcommand{\exmpl}{\end{xmpl}}
\newcommand{\ermrk}{\end{rmrk}}

\newcommand{\bbc}{\mathbb{C}}
\newcommand{\bbz}{\mathbb{Z}}
\newcommand{\bbm}{\mathbb{M}}
\newcommand{\bbn}{\mathbb{N}}

\newcommand{\bbt}{\mathbb{T}}

\newcommand{\cla}{\mathcal{A}}

\newcommand{\clh}{\mathcal{H}}

\newcommand{\clk}{\mathcal{K}}
\newcommand{\cll}{\mathcal{L}}

\newcommand{\prf}{\noindent{\it Proof\/}: }
\newcommand{\seq}{\subseteq}
\newcommand{\one}{{1\!\!1}}
\newcommand{\nn}{\nonumber}

\newcommand{\id}{\mbox{id}}

\def \qed { \mbox{}\hfill 
$\Box$\vspace{1ex}}

\newcommand{\sgn}{\mbox{sign\,}}
\newcommand{\half}{\frac{1}{2}}

\begin{document}
\author{{\sc Partha Sarathi Chakraborty} and
{\sc Arupkumar Pal}}
\title{An invariant for homogeneous spaces of compact quantum groups}
\maketitle

\begin{abstract}
The central notion in Connes' formulation of non commutative geometry
is that of a spectral triple.
Given a homogeneous space of a compact quantum group, restricting our attention
to all spectral triples that are `well behaved' with respect to
the group action, we construct a certain dimensional invariant.
In particular, taking the (quantum) group itself as the homogeneous space,
this gives an invariant for a compact quantum group.
Computations of this invariant in several cases, including all
type A quantum groups, are given.
\end{abstract}
\maketitle

\section{Introduction}
 Geometry can be broadly interpreted as the study of cycles and their
intersection properties in some suitable homology theory.
Noncommutative Geometry is no exception. In Alain Connes'
interpretation, noncommutative geometry is the study of spectral triples
or unbounded Kasparov-modules  with finer properties (\cite{con-1980a, con-1994a}). Often
these finer properties encode information about metric, dimension etc.
In fact for an unbounded K-cycle to encode useful information, it is not
always necessary that the cycle should be homologically nontrivial. A
cycle may be homologically trivial but still it may contain metric or
dimensional information. One prime example is the Laplacian on odd
dimensional manifolds. Study of cycles not necessarily nontrivial is 
not new.  Apart from Connes (see for example \cite{con-1989a}), it also 
includes Voiculescu's work (\cite{voi-1979a}, \cite{voi-1981a}) 
on norm ideal perturbations or Rieffel's work (\cite{rie-2004a}) on 
extending the notion of metric spaces. Voiculescu
answers the question of existence of bounded K-cycles in a given
representation of a $C^*$-algebra. He does not comment on the
nontriviality of the K-cycle as a K-homology class. Rieffel uses
spectral triples to produce compact quantum metric spaces but in his
construction nontriviality do not play any major role.
Here we take a similar approach. We utilise the notion of spectral
triples to produce dimensional invariants for ergodic $C^*$-dynamical
systems. However, we should emphasize a crucial distinction between the
above cited works and the present one. In Rieffel's work, Dirac operators are used
as a source to produce compact quantum metric spaces; but there are situations
(\cite{cha-2005a}) where one produces compact quantum metric spaces
even without using spectral triples. Whereas in our case, the concept of spectral
triple is used in an essential manner, not as a source of examples.

Origin of the present paper lies in the search for non trivial spectral triples
for the quantum $SU(2)$ (\cite{cha-pal-2003a}) and quantum spheres (\cite{cha-pal-2008a}).
Instead of KK-theoretic machinery, our method tries to characterize equivariant 
spectral triples. In the process, one observes that even without the condition
of nontriviality of the corresponding K-cycle, in certain cases there
are canonical spectral triples encoding essential information about the space.
That leads to the present invariant.

Let us recall the classical fact that if $M$ is a $d$-dimensional
compact Riemannian manifold and we have an elliptic operator $D$ of
order 1 then the $n$-th singular value of $|D|$ grows like $n^{1/d}$.
Using this Connes has indicated how to define the dimension of a
spectral triple. This observation can be utilized to obtain invariants
for $C^*$-algebras provided one could associate natural spectral
triples with them. 
But unlike the classical case it
is difficult to define natural spectral triples for $C^*$-algebras. As
one tries to answer this question one faces two main  dfficulties: one, it
is difficult to get hold of a canonical representation other than the
GNS representation, which normally is too big; two:  
once a `natural' representation has been identified, it is 
impossible to classify all spectral triples so as to be able to extract
any meaningful quantity out of them.
Another important point pertaining to the last problem above is
that when considering a spectral triple, one should look at a spectral triple for
some dense *-subalgebra of the $C^*$-algebra, but in general there is no
canonical choice of a dense *-subalgebra; but properties of the behaviour
of the spectral triple is sensitive to this choice.

If we restrict our attention to
homogeneous spaces of compact quantum groups, or which 
is the same thing, to ergodic $C^*$-dynamical systems,
then these problems can be resolved to a large extent. 
Ergodicity gives us a unique invariant state, so  passing to its GNS 
representation gives us a way of fixing a canonical representation. 
Once that is done, we can again use the group action to identify an appropriate
dense *-subalgebra.
As we shall see, in many situations, it is then possible to 
investigate a subclass of spectral triples for this dense *-algebra
that are `well-behaved' with respect to the structure at hand;
and one can study them in  enough details in order to be able to come up with an invariant.

In the next section, we make all these notions precise and define the
invariant. In the remaining sections, we compute it 
in several cases. In particular, in section~5, this is done for all
type A quantum groups. When one introduces an invariant for a class
of objects, one of the major issues one must look at is its
computability. Thus section~5 can be looked upon as the heart of
the present paper.

\section{The Invariant}
Let us start by recalling the notion of a spectral triple.

  \bdfn
Let $\cla$ be an associative unital  *-algebra.
A \textbf{spectral triple} for $\cla$ is a triple $(\clh,\pi,D)$ where
\begin{enumerate}
 \item
$\clh$ is a (complex separable) Hilbert space,
\item
$\pi:\cla\longrightarrow \cll(\clh)$ is a *-representation (usually assumed faithful),
\item
$D$ is a self-adjoint operator with compact resolvent such that
$[D,\pi(a)]\in\cll(\clh)$ for all $a\in \cla$.
\end{enumerate}
One often writes $(\clh,\cla,D)$ in place of $(\clh,\pi,D)$
if the representation $\pi$ is clear from the context.
We say that a spectral triple is $p$-summable if ${|I+D^2|}^{-p/2}$
belongs to the ideal $\cll^1$ of Trace-class operators.
\edfn
\brmrk
Since $D$ has compact resolvent, the above is equivalent to
saying that $|D|^{-p}$ is Trace-class on the complement
of its kernel.
In some of our earlier papers (\cite{cha-pal-2003a, cha-pal-2008a}), the phrase
`$p$-summable' was used to mean that $|I+D^2|^{-p/2}$ is in
the ideal $\cll^{(1,\infty)}$ of Dixmier traceable operators, which is
slightly different from the definition that we use here.
\ermrk 
 
It is a well known fact \cite{gilkey-1995a} that if
$M$ is a compact Riemannian manifold  of dimension $d$ and $D$ is a
elliptic differential operator of order $r$ then $\mu_n(|D|)$, the
$n$-th smallest singularvalue of $D$ behaves like $n^{\frac{r}{d}}$.
Therefore one has the equality 
\[
d=\inf\{ \delta:
\mbox{Tr}{|D|}^{-\frac{\delta}{r}} < \infty \}.
\]
Motivated by this Alain Connes defined the dimension of a spectral
triple ${\mathcal E}=(\cla,\clh,D)$ as 
\[
\mbox{dim\,}{\mathcal
E}=\inf\{\delta: \mbox{Tr}{|D|}^{-\delta} < \infty\}.
\]
We want to utilize this concept to define a dimensional invariant for a
$C^*$-algebra. There are two obstacles, firstly there is no canonical
representation associated with a $C^*$-algebra and secondly even after
fixing a representation there is no canonical Dirac operator. One simple
strategy to overcome the first problem is to consider $C^*$-dynamical
systems $(A,G,\tau)$ with an ergodic action of a compact quantum group
$G$. Then by passing to the GNS representation of the invariant state we
get a canonical representation.
 To tackle the second problem we exploit the notion of equivariance.
 
 Let us now recall  various terminologies used here
 
\bdfn \label{df:cqg}(\cite{wor-1998a})
	A compact quantum group $G$ consists of a unital $C^*$-algebra 
$C(G)$ along with a unital homomorphism $\Delta: C(G) \rightarrow C(G)
\otimes C(G)$ such that 
\begin{enumerate}
\item 
$(\Delta \otimes \id)\Delta=(\id \otimes
\Delta)\Delta$, 
\item
 both $\{(a\otimes I)\Delta(b): a, b\in C(G)\}$ and 
$\{(I \otimes a)\Delta(b): a,b\in C(G)\}$ are total  in $C(G)$.
\end{enumerate}
\edfn
It is known that any compact quantum group admits a unique invariant state $h$,
called the Haar state. The invariance property that the Haar state has is the following:
\[
(h\otimes id)\Delta (a)=h(a)I = (id\otimes h)\Delta (a) \quad\mbox{for all } a\in C(G).
\]

\bdfn\label{df:cqgaction} (\cite{pod-1995a})
We say that  a compact quantum group $G$ acts on a $C^*$-algebra $A$ if there
is a homomorphism $\tau: A \rightarrow A \otimes C(G)$ such that 
\begin{enumerate}
\item 
$(\tau \otimes \id)\tau=(\id \otimes
\Delta)\tau$, 
\item
$\{(I \otimes b)\Delta(a): a\in A, b\in C(G)\}$ are total  in $A\otimes C(G)$.
\end{enumerate}
The $C^*$-algebra $A$ is called a \textbf{homogeneous space} of $G$ if
the fixed point subalgebra $\{a\in A: \tau(a)=a\otimes I\}$ is $\mathbb{C}I$.
In such a case, the action $\tau$ is said to be \textbf{ergodic} and we call $(A, G,\tau)$
an ergodic $C^*$-dynamical system.
\edfn
Recall (Proposition~1.9, \cite{pod-1995a}) that if $A$ happens to be
a quotient space of $G$ (i.e.\ $A$ is isomorphic to $C(G\verb1\1H)$ for some
closed quantum subgroup $H$ of $G$ and the $G$-action on $A$ is
equivalent to the natural $G$-action on $G\verb1\1H$), then $A$ is a homogeneous 
space of $G$ and the action is ergodic.

A covariant representation $(\pi,u)$
of a $C^*$-dynamical system $(A,G,\tau)$ consists of
a unital *-representation $\pi:A\rightarrow\cll(\clh)$,
a unitary representation $u$ of $G$ on $\clh$, i.e.\
a unitary element of the multiplier algebra $M(\clk(\clh)\otimes C(G))$
such that they obey the condition
$(\pi\otimes\id)\tau(a)=u(\pi(a)\otimes I)u^*$ for all $a\in A$.
\bdfn
Suppose $(A, G,\tau)$ is a $C^*$-dynamical system.
An operator $D$ acting on a Hilbert space $\clh$
is said to be \textbf{equivariant} with respect to a covariant
representation $(\pi,u)$ of the system if
$D\otimes I$ commutes with $u$. If $(\pi,u)$ is a covariant representation
of $(A,G,\tau)$ on a Hilbert space $\clh$ and $(\clh,\pi,D)$ is a spectral
triple for a dense *-subalgebra $\cla$ of $A$, then we say that 
$(\clh,\pi,D)$ is equivariant with respect to
$(\pi,u)$ if the operator $D$ is equivariant with respect to $(\pi,u)$.
\edfn

  A homogeneous space $A$ for $G$ admits an invariant state
  $\rho$ that satisfies
  \[
  (\rho\otimes id)\tau(a)=\rho(a)I,\quad a\in A.
  \]
  This invariant  state $\rho$ is unique and is related to the 
  Haar state $h$ on $G$ through the
  equality 
  \[
  (id\otimes h)\tau (a)= \rho(a)I, \quad a\in A.
  \]

Given an ergodic $C^*$-dynamical system $(A,G,\tau)$ with unique
invariant state $\rho$, denote by $(\clh_\rho,\pi_\rho,\eta_\rho)$ the GNS representation
associated with the state $\rho$, i.e.\ $\clh_\rho$ is a Hilbert space,
$\eta_\rho:A\rightarrow\clh_\rho$ is linear with
$\eta_\rho(A)$ dense in $\clh_\rho$ and $\langle \eta_\rho(a),\eta_\rho(b)\rangle=\rho(a^*b)$;
and $\pi_\rho:A\rightarrow\cll(\clh_\rho)$ is the *-representation of $A$ on $\clh_\rho$
given by $\pi_\rho(a)\eta_\rho(b)=\eta_\rho(ab)$. The action $\tau$ induces a unitary representation
$u_\tau$ of $G$ on $\clh_\rho$ that makes the pair $(\pi_\rho,u_\tau)$ a covariant
representation of the system $(A,G,\tau)$. We fix this covariant representation.
Let $A(G)$ denote the dense *-subalgebra of $C(G)$ generated by the matrix entries
of irreducible unitary representations of $G$, and let
\[
\cla=\{a\in A: \tau(a)\in A\otimes_{alg} A(G)\}.
\]
By the results in \cite{pod-1995a}, $\cla$ is a dense *-subalgebra 
of $A$.
Now  consider the class ${\mathcal E}$ of  spectral triples for $\cla$ equivariant
with respect to the covariant representation $(\pi_\rho,u_\tau)$. 
We define the spectral dimension
of the system $(A,G,\tau)$ to be the quantity 
\[
\inf\{p>0: \exists\,D \text{ such that }(\clh_\rho,\pi_\rho,D)\in\mathcal{E} 
   \text{ and $D$ is $p$-summable}\}.
\]
We will denote this number by ${\mathcal Sdim}(A,G,\tau)$.
Here we have taken infimum because if we have a spectral triple 
$\mathcal{D}=(\clh,\pi,D)$ that is $p$-summable,
then  ${\mathcal D}_\alpha=(\clh,\pi, D/{|D|}^\alpha), 0 <
\alpha < 1$ is a spectral triple that is  $p/(1-\alpha)$-summable.
(see page~459, \cite{g-v-f-2001a}; \cite{con-mos-1986a}). 

The first instance of computation of this invariant
can be traced back to Connes in \cite{con-1989a}, where
he proved that, in the terminology of the present paper,
if $\Gamma$ is a discrete group containing the free group
on two generators, then the Pontryagin dual  
$\hat{\Gamma}$, which is a compact quantum group, has spectral dimension $\infty$.
In the rest of the paper we will focus on computation of this
invariant in several other cases, emphasis being on
examples where the number is finite.

\section{A Commutative Example }
In this section we will  calculate the spectral dimension 
of $SU(2)$ under its own natural action.

Recall from \cite{cha-pal-2003a} that the representation of $C(SU(2))$ on $L_2(SU(2))$ is given by
\bean
\alpha: e^{(n)}_{ij}  &\mapsto& a_+(n,i,j) e^{(n+\half)}_{i-\half ,j-\half } 
       + a_-(n,i,j)  e^{(n-\half )}_{i-\half ,j-\half },\label{calpha}\\
\beta:e^{(n)}_{ij}  &\mapsto& b_+(n,i,j)  e^{(n+\half )}_{i+\half ,j-\half } 
       + b_-(n,i,j)  e^{(n-\half )}_{i+\half ,j-\half },\label{cbeta}
\eean
where $n \in \half\bbn $, $i,j \in \{ -n,-n+1,\cdots,n\}$  and
\bea
a_+(n,i,j)  & =& \Bigl(\frac{(n-j+1)(n-i+1)}{(2n+1)(2n+2)}
                          \Bigr)^\half,\label{caplus}\\
a_-(n,i,j)&=&\Bigl(\frac{(n+j)(n+i)}{2n(2n+1)}\Bigr)^\half,\label{caminus}
\\
b_+(n,i,j)&=& -
\Bigl(\frac{(n-j+1)(n+i+1)}{(2n+1)(2n+2)}\Bigr)^\half,\label{cbplus}\\
b_-(n,i,j)&=&\Bigl(\frac{(n+j)(n-i)}{2n(2n+1)}\Bigr)^\half,\label{cbminus}
\eea

Observe that the representation of (the complexification of)
$\mathfrak{su}(2)$ 
on $L_2(SU(2))$ is given by
\bean
h e^{(n)}_{ij} &=& (n-2j)e^{(n)}_{ij},\\
e e^{(n)}_{ij} &=& j(n-2j+1)e^{(n)}_{i,j-1},\\
f e^{(n)}_{ij} &=& e^{(n)}_{i,j+1},
\eean
where $h$, $e$ and $f$ obey
\[
[h,e]=2e,\quad [h,f]=-2f,\quad [e,f]=h.
\]
Therefore,  any equivariant self-adjoint operator with discrete
spectrum must be of the form
\be\label{genericeq}
D: e^{(n)}_{ij}\mapsto d(n,i)e^{(n)}_{ij}.
\ee
Commutators of this operator with  $\alpha$ and $\beta$
are  given by
\bea 
[D,\alpha]e^{(n)}_{ij} &=& a_+(n,i,j)(d(n+\half,i-\half)-d(n,i))
                                                      e^{(n+\half)}_{i-\half
,j-\half}\nn\\
&&{}     +a_{-}(n,i,j)(d(n-\half,i-\half)-d(n,i))e^{(n-\half)}_{i-\half,j-\half},
\label{eq:bdd1}\\
{}[D,\beta]e^{(n)}_{ij} &=& b_+(n,i,j)(d(n+\half,i+\half)-d(n,i))
                                                  e^{(n+\half)}_{i+\half,j-\half}\nn\\
      &&{}   
+b_{-}(n,i,j)(d(n-\half,i+\half)-d(n,i))e^{(n-\half)}_{i+\half,j-\half}.\label{eq:bdd2}
\eea
where $a_{\pm}$ and $b_{\pm}$ are now given by equations
(\ref{caplus})--(\ref{cbminus}).

\blmma[\cite{cha-pal-2003a}]\label{sign-sum}
Suppose $D$ is an operator on $L_2(SU(2))$ given by (\ref{genericeq})
and having bounded commutators with $\alpha$ and $\beta$.
Then $D$ can not be $p$-summable for $p\le 3$.
\elmma
\prf
Conditions for boundedness of the commutators give us
\bea
|d(n+\half,i+\half)-d(n,i)| &=&
O\Bigl(\left(\frac{2n+2}{n+i+1}\right)^\half\Bigr),\label{cbdd1}\\
|d(n+\half,i-\half)-d(n,i)| &=&
O\Bigl(\left(\frac{2n+2}{n-i+1}\right)^\half\Bigr).\label{cbdd2}
\eea
Observe from (\ref{cbdd1}) and (\ref{cbdd2}) that if we restrict
ourselves to 
the region $i\geq 0$, then
\be\label{cbdd3}
|d(n+\half,i+\half)-d(n,i)| = O(1),
\ee
and if we restrict to $i\leq 0$, then 
\be\label{cbdd4}
|d(n+\half,i-\half)-d(n,i)| = O(1).
\ee
Also, it is not too difficult to see that
\be\label{cbdd5}
|d(n+1,0)-d(n,0)| = O(1).
\ee
Suppose $C>0$ is a constant that works for (\ref{cbdd1})--(\ref{cbdd5}).
Then
\be\label{bndgenericdq1}
|d(n,i)| < 2C n
\ee
Therefore, for $\alpha \leq 3$ 
\bean
   \mbox{Tr} |D|^{-\alpha} > \sum_{n \ge 0} {(2n+1)}^2 {(2Cn)}^{-\alpha}
= \infty
\eean
\qed
\bthm
Spectral dimension of $SU(2)$ is 3.
\ethm
\prf The previous lemma shows that ${\mathcal Sdim}(C(SU(2)))\ge3$. Note
that $D$ given by $d(n,i)=n$ will give rise to a spectral triple that is $p$-summable
for all $p>3$. Hence the result.\qed


\section{$\bbt^n$ action on noncommutative torus}
Our next example deals with probably the most well known ergodic
$C^*$-dynamical system, namely that of $\bbt^n$ acting ergodically on
the noncommutative torus.  By a result of Milnes and Walters
(\cite{mil-wal-2000a}), this describes all the primitive $C^*$-algebras with free
and ergodic action of $\bbt^n$. Recall (\cite{rie-1990a}) that the
noncommutative $n$-torus is the universal $C^*$-algebra generated by
$n$-unitaries $U_1,\cdots, U_n$ satisfying the commutation relation
$U_jU_k=exp(2\pi i \theta_{jk})U_kU_j$, where $\Theta=(\!(\theta_{jk})\!)$
is a skew symmetric matrix with real entries.
This $C^*$-algebra is referred as the noncommutative torus and denoted
by $A_\Theta$. If $\Theta$ has sufficient irrationality (i.e.\ 
$p^t\Theta q\in\bbz$ for all $q\in{\mathbb Z}^n$ implies $ p=0$), 
then $A_\Theta$ is simple.
Non-commutative torus admits an action of $\bbt^n$. To specify the
action it is enough to prescribe them on the generators by
$\alpha_{z}(U_j)=z_j U_j, 1 \le j \le n, z=(z_1,\cdots,z_n) \in \bbt^n$.
This action is ergodic. Therefore the dynamical system $(A_\Theta,
\bbt^n,\alpha)$ satisfies our hypothesis and we can ask what is the
spectral dimension of this system.
\bppsn
The spectral dimension of $(A_\Theta, \bbt^n,\alpha)$ is n.
\eppsn 
\prf 
The unique invariant state is specified by
$$\tau(a)=\int_{\bbt^n} \alpha_z(a) dz.$$
The algebra of smooth elements $A_\Theta^\infty=\{a \in A_\Theta: z
\mapsto \alpha_z(a) \mbox{ is smooth}  \}$ is dense in $A_\Theta$ and
can  be described as $\{ a: a=\sum a_{k_1,k_2,\cdots,k_n}
U_1^{k_1}\cdots U_n^{k_n}, a_{\underline{k}} \in S({\mathbb Z}^n) \}$.
On $A_\Theta^\infty$ the invariant state is specified by $\tau(\sum
a_{k_1,k_2,\cdots,k_n} U_1^{k_1}\cdots U_n^{k_n})=a_{0,\cdots,0}$. This
state is tracial and the associated GNS space can be identified as
$l_2({\mathbb Z}^n)$. If we denote by $\{ e_{k_1,\cdots,k_n} :
k_1,\cdots,k_n \in {\mathbb Z} \}$ the canonical basis elements then the
GNS representation is given by 
\[
U_j(e_{k_1,\cdots,k_n})=\exp \left(2 \pi i
\sum_{r=1}^{j-1} \theta_{jr}k_r\right)e_{\underline{k}+\epsilon_j},
\]
where
$\epsilon_j=(0,\cdots,1,\cdots,0)$, the $1$ is at the j-th position. The
$e_{\underline{k}}$ are the spectral subspaces and hence an equivariant
$D$ must be of the form 
\[
D: e_{\underline{k}} \mapsto
d(\underline{k})e_{\underline{k}}.
\]
The boundedness of the commutator condition is equivalent to
\be\label{bddnct}
|d(\underline{k}+\epsilon_j)-d(\underline{k})| < C, \quad\forall
\;\underline{k} \in {\mathbb Z}^n, \quad 1 \le j \le n
\ee
for some constant $C$. Therefore
$$|d(\underline{k})| < C |\underline{k}|,$$
where $|\underline{k}|=k_1+\cdots + k_n$. Then 
\bean
\mbox{Tr} {|D|}^{-\alpha}&=& \sum_{k_1,\cdots,k_n \in {\mathbb Z}}
C^{-\alpha} {|k|}^{-\alpha}\\
& >& {C'}^{-\alpha}\sum_{k_1,\cdots,k_n
=2}^\infty \left(\sum_{j=1}^n k_j^2\right)^{-\frac{\alpha}{2}} \\
&\ge& \int_{1}^{\infty}\cdots \int_1^\infty \left(\sum_{j=1}^n
x_j^2\right)^{-\frac{\alpha}{2}}dx_1\cdots dx_n\\
&>& \int_{\sqrt{n}}^\infty \int_{S^{n-1}} r^{-\alpha} r^{n-1} dr
d\sigma,
\eean
 where $d\sigma$ is the normalized surface measure on the sphere.
Therefore $\mbox{Tr} {|D|}^{-\alpha} $ is finite only if
$\int_{1}^\infty r^{-\alpha+n-1} dr < \infty$, that is $\alpha> n$.

On the other hand by taking $d(\underline{k})=|\underline{k}|$ we obtain
a $D$ such that $\mbox{Tr} |D|^{-\alpha} < \infty$ for all $\alpha > n$.
Hence spectral dimension of $A_\Theta$ is $n$. \qed

\section{$SU_q(\ell+1)$ action on itself}
In this section and the next, we will discuss two $C^*$-dynamical systems,
both involving the quantum group $SU_q(\ell+1)$. In one case, the
$C^*$-algebra under consideration will be the $C^*$-algebra $C(SU_q(\ell+1))$
of `continuous functions' associated with the quantum group $SU_q(\ell+1)$
itself, and in the other case, it will be the $C^*$-algebra
$C(S_q^{2\ell+1})$ of `continuous functions' on the odd dimensional
spheres $S_q^{2\ell+1}$. In both sections, $q$ is assumed to lie in
the open interval $(0,1)$. 

We will take $G$ to be the quantum group $SU_q(\ell+1)$
and $A$ to be the $C^*$-algebra  $C(G)$. Action of $G$ on $A$ will
just be the comultiplication of $G$. Haar state of $G$ is the invariant state for this
action and therefore the relevant covariant representation of this
$C^*$-dynamical system is given by the triple
$(L_2(G), \pi,u)$ where $L_2(G)$ is the GNS space of the Haar state on $A=C(G)$,
$\pi$ is the representation of $A$ on $L_2(G)$ by left multiplication, and
$u$ is the right regular representation. 

The $C^*$-algebra $C(SU_q(\ell+1))$ is the universal $C^*$-algebra
generated by  $\{ u_{ij} : i,j=1, \cdots ,\ell+1 \} $ 
obeying the relations (see \cite{wor-1988a}):
\[ 
\sum_k u_{ki}^* u_{kj}= \delta_{ij} I,\quad
\sum_k u_{ik} u_{jk}^*= \delta_{ij} I
\] 
\[
\sum_{\mbox{\scriptsize{$k_i$'s distinct}}} 
\hspace{-1em}(-q)^{I(k_1,k_2, \cdots, k_{\ell+1})} 
   u_{j_1k_1}
\cdots u_{j_{\ell+1} k_{\ell+1}} = \begin{cases}
      (-q)^{I( j_1,j_2, \cdots, j_{\ell+1})} & \mbox{$j_i$'s distinct}\cr
0 & \mbox{otherwise}
             \end{cases}
\]
where $I( k_1,k_2, \cdots, k_{\ell+1})$ is the number of inversions
in  $( k_1,k_2, \cdots, k_{\ell+1})$.
The group laws are given by the folowing maps:
\bean
\Delta(u_{ij})&=&\sum_k u_{ik}\otimes u_{kj}\qquad\mbox{ (Comultiplication)}\\
 S(u_{ij})&=&u_{ji}^*\qquad\mbox{ (Antipode)}\\
  \epsilon(u_{ij})&=&\delta_{ij}\qquad\mbox{ (Counit)}
\eean

For preliminaries on irreducible unitary representations of the group
$SU_q(\ell+1)$, we refer the reader to section~7.3, \cite{kli-sch-1997a} and
section~4, \cite{cha-pal-2008a}.
For related computations that will be very crucial for us, we refer the
reader to \cite{cha-pal-2008a} and \cite{pal-sun-2010a}.
We will stick to the same notations as in those two papers.
For convenience, let us summarize the main points.
Let $\Lambda$ denote the set of all Young tableaux 
\[
\{\lambda\equiv(\lambda_1,\ldots,\lambda_{\ell+1}\}:\lambda_i\in\bbn, \lambda_1\geq\lambda_2\geq\ldots\lambda_{\ell+1}=0\}
\]
and for $\lambda\in\Lambda$, let $\Gamma_\lambda$ denote the set of all 
Gelfand-Tsetlin tableaux (will be referred to as GT
tableaux) with top row equal to $\lambda$.

\begin{enumerate}
\item 
Irreducible unitary representations of $G$ are indexed by $\Lambda$. We will denote them by $u^\lambda$.
\item 
dimension of the Hilbert space on which $u^\lambda$ acts is $|\Gamma_\lambda|$,
\item 
Fixing an orthonormal basis $\{e_\bldr: \bldr\in \Gamma_\lambda\}$ of this Hilbert space,
one gets the matrix entries $u^\lambda_{\bldr \blds}$, $\bldr,\blds\in \Gamma_\lambda$.
\item 
$\{u^\lambda_{\bldr \blds}: \lambda\in \Lambda, \bldr,\blds\in \Gamma_\lambda\}$ generate the $C^*$-algebra $A=C(G)$,
\item 
denote by $e^\lambda_{\bldr \blds}$ the normalized $u^\lambda_{\bldr \blds}$'s, i.e.\
$e^\lambda_{\bldr \blds}=\|u^\lambda_{\bldr \blds}\|^{-1}u^\lambda_{\bldr \blds}$. Then
$\{e^\lambda_{\bldr \blds}: \lambda\in \Lambda, \bldr,\blds\in \Gamma_\lambda\}$ is a complete orthonormal
basis for $L_2(G)$.
\item 
For the Young tableaux
$\one:=(1,0,\ldots,0)$. We will omit the symbol $\one$
and just write $u$ in order to denote $u^\one$.
Notice that any GT tableaux $\bldr$ with first row $\one$
must be, for some $i\in\{1,2,\ldots,\ell+1\}$, of the form $(r_{ab})$, where
\[
r_{ab}=\begin{cases}1 & \mbox{if $1\leq a\leq i$ and $b=1$},\cr
       0   &\mbox{otherwise.}\end{cases}
\]
Thus such a GT tableaux is uniquely determined by the integer $i$.
We will write just $i$ for this GT tableaux $\bldr$.
Thus a typical matrix entry of $u^\one$ will be
written simply as $u_{ij}$.
\item
Let $\bbm_k:=\{M=(m_1,m_2,\ldots,m_k)\in\bbn^k: 1\leq m_j\leq \ell+2-j, 1\leq j\leq k\}$.
For $M\in\bbm_k$ and a GT tableaux $\bldr$, denote by $M(\bldr)$ the tableaux $\blds$ defined by
\be\label{eq:movenotation}
s_{ij}=\begin{cases}
         r_{ij}+1 & \mbox{ if } j=m_i, 1\leq i\leq k,\cr
         r_{ij} & \mbox{ otherwise},
         \end{cases}
\ee
and let $\sgn(M)$ denote
the product $\prod_{a=1}^{k-1}\sgn(m_{a+1}-m_a)$ (here $\sgn(p):=2\chi_\bbn(p)-1$).
Then one has
\be\label{eq:left_mult}
\pi(u_{ij})e_{\bldr\blds}
= \sum_{\genfrac{}{}{0pt}{}{M\in\bbm_i, N\in\bbm_j}{m_1=n_1}}
  \sgn(M)\sgn(N)q^{C(i,\bldr,M)+C(j,\blds,N)+A(M)+K(M)+B(N)}(1+o(q))
   e_{M(\bldr)N(\blds)},
\ee
where
\bea
C(i,\bldr,M) &= &\sum_{a=1}^{i-1}\left(
   \sum_{b=m_a\wedge m_{a+1}}^{m_a\vee m_{a+1} -1}H_{ab}(\bldr)
   +2 \sum_{m_{a+1} < b < m_a}V_{ab}(\bldr)\right)
+\sum_{m_i \leq b < \ell+2-i}H_{ib}(\bldr),\nonumber\\
&&\label{eq:cgc5}\\
 A(M) &=& \sum_{j=1}^{i-1}\left|m_j - m_{j+1}\right|-\#\{1\leq j\leq
i-1:m_j> m_{j+1}\},\label{eq:A}\\
K(M) &=& \ell+2-i-m_i,\label{eq:K}\\
B(M) &=& A(M)+m_1-m_i. \label{eq:B}
\eea
\end{enumerate}

Let us next derive the general form of the operator $D$ from
the equivariance condition.
The Hilbert space $\clh=L_2(G)$ decomposes as a direct sum
$\oplus_\lambda\clh_\lambda$
where $\clh_\lambda=\mbox{span}\{e^\lambda_{\bldr,\blds}:\bldr,\blds\in\Gamma_\lambda\}$
and the restriction of $u$ to each $\clh_\lambda$ is equivalent to a
direct sum of $\dim\lambda$ copies of $u^\lambda$ and the
operator $D$ respects this decomposition.
Looking at the restriction of $D$ to each $\clh_\lambda$ and
using the fact that it commutes with $u$, it follows that
the restriction of $D$ to $\clh_\lambda$ is of the form
$\oplus_\mu d_{\lambda\mu}P_{\lambda\mu}$
where $u$ commutes with each $P_{\lambda\mu}$ and
the restriction of $u$ to each $P_{\lambda\mu}$ is equivalent to
$u^\lambda$.
Write the restriction of $D$ to $\clh_\lambda$ in the form
$\sum_\mu d_{\lambda\mu}P_{\lambda\mu}$ where
the $d_{\lambda\mu}$'s are distinct for distinct $\mu$'s.
Write $\clh_{\lambda\mu}=P_{\lambda\mu}\clh_\lambda$.
\bppsn\label{pr:equivariantD-1}
Choose and fix a $\mu$. 
Then
\begin{enumerate}
\item 
If $\sum_{\blds,\bldt}c(\blds,\bldt)e^\lambda_{\blds,\bldt}\in \clh_{\lambda\mu}$,
 then $\sum_{\blds}c(\blds,\bldt)e^\lambda_{\blds,\bldt}\in \clh_{\lambda\mu}$ for all $\bldt\in\Gamma_\lambda$.
 \item
 If $\sum_{\blds}c(\blds)e^\lambda_{\blds,\bldt_0}\in \clh_{\lambda\mu}$ for some $\bldt_0\in\Gamma_\lambda$,
 then $\sum_{\blds}c(\blds)e^\lambda_{\blds,\bldt}\in \clh_{\lambda\mu}$ for all $\bldt\in\Gamma_\lambda$.
\end{enumerate}
\eppsn
\prf
1. Choose and fix an $\bldt_0\in\Gamma_\lambda$.
Take a linear functional on $\clh_\lambda\seq C(G)$ that takes
the value 1 at $e^\lambda_{\bldt_0\bldt_0}$ and vanishes at all other
$e^\lambda_{\blds\bldt}$'s. Extend it to a bounded linear functional $\rho$ on $C(G)$.
Write $u_\rho:=(\id\otimes\rho)u$.
Then 
\[
u_\rho\left(\sum_{\blds,\bldt}c(\blds,\bldt)e^\lambda_{\blds,\bldt}\right)=
     \sum_{\blds}c(\blds,\bldt_0)e^\lambda_{\blds,\bldt_0}.
\]
Since $D$ commutes with $u_\rho$ and the $d_{\lambda\mu}$'s are distinct,
it follows that $\sum_{\blds}c(\blds,\bldt_0)e^\lambda_{\blds,\bldt_0}\in \clh_{\lambda\mu}$.

2. Choose and fix $\bldt_1\in\Gamma_\lambda$. 
In this case, take a linear functional $\rho$ on $C(G)$ such that
\[
\rho(e^\lambda_{\bldt_1\bldt_0})=1, \quad \rho(e^\lambda_{\blds\bldt})=0 \quad\mbox{for all }(\blds,\bldt)\neq (\bldt_1,\bldt_0).
\]
Then
$u_\rho\left(\sum_{\blds}c(\blds)e^\lambda_{\blds,\bldt_0}\right)=\sum_{\blds}c(\blds)e^\lambda_{\blds,\bldt_1}\in\clh_{\lambda\mu}$.
\qed

\bcrlre\label{cor:equivariantD-2}
There exists a unitary $(\!(c^\bldr_\blds)\!)_{\bldr,\blds}\in GL(|\Gamma_\lambda|,\bbc)$
and a partition $\{F_\mu:\mu\}$ of $\Gamma_\lambda$
such that
\be\label{eq:equivD-1}
\clh_{\lambda\mu}=\mbox{span}\left\{\sum_\blds c^\bldr_\blds e^\lambda_{\blds\bldt}: \bldr\in F_\mu, \bldt\in\Gamma_\lambda\right\}.
\ee
\ecrlre
\prf
It follows from the above proposition that there is a subset $E_\mu$ of
$\bbc^{|\Gamma_\lambda|}$ such that
\[
\clh_{\lambda\mu}=\mbox{span}\left\{\sum_\blds c(\blds)e^\lambda_{\blds\bldt}: (c(\blds))_\blds\in E_\mu, \bldt\in\Gamma_\lambda\right\}.
\]
One can now further assume, without loss in generality, that
$E_\mu$ is an independent set of unit vectors in $\bbc^{|\Gamma_\lambda|}$. 
Orthogonalizing the vectors $(c(\blds))_\blds$ next, we obtain
a set of orthonormal vectors $E_\mu$ in $\bbc^{|\Gamma_\lambda|}$
such that
\[
\clh_{\lambda\mu}=\mbox{span}\left\{\sum_\blds c(\blds)e^\lambda_{\blds\bldt}: (c(\blds))_\blds\in E_\mu, \bldt\in\Gamma_\lambda\right\}.
\]
Since the spaces $\clh_{\lambda\mu}$ are orthogonal for different $\mu$,
$\clh_\lambda=\oplus_\mu\clh_{\lambda\mu}$ and $\dim\clh_{\lambda\mu}=|\Gamma_\lambda|\cdot|E\mu|$,
it follows that the collection of vectors $(c(\blds))_\blds \in\cup_\mu E_\mu$ form
a complete orthormal set of vectors in $\bbc^{|\Gamma_\lambda|}$.
The matrix formed by taking the elements of $\cup_\mu E_\mu$ as rows
now give us the required conclusion, with $F_\mu$ being the rows
corresponding to $(c(\blds))_\blds\in E_\mu$.
\qed

\bppsn\label{pr:equivariantD-3}
If we make a change of basis in the Hilbert space on which the irreducible
$u^\lambda$ acts using the matrix $(\!(c^\bldr_\blds)\!)$, then with respect to
the new matrix entries $u^\lambda_{\bldr\blds}$ and 
$e^{\lambda}_{\bldr\blds}=\|u^\lambda_{\bldr\blds}\|^{-1}|u^\lambda_{\bldr\blds}$,
one has
\[
\clh_{\lambda\mu}=\mbox{span}\left\{e^\lambda_{\bldr\blds}:\bldr\in F_\mu,\blds\in\Gamma_\lambda\right\}.
\]
\eppsn
\prf
This follows from the observation that
\[
\mbox{span}\left\{\sum_\blds c^\bldr_\blds e^\lambda_{\blds\bldt}: \bldt\in\Gamma_\lambda\right\}
  =\mbox{span}\left\{\sum_{\blds,\bldt} c^\bldr_\blds e^\lambda_{\blds\bldt} \overline{c^\bldz_\bldt}: \bldz\in\Gamma_\lambda\right\},
\]
and $\sum_{\blds,\bldt} c^\bldr_\blds e^\lambda_{\blds\bldt} \overline{c^\bldz_\bldt}$
are the matrix entries with respect to the new basis
$f_\bldr:=\sum_\blds c^\bldr_\blds e_\blds$.
\qed

Thus we can assume that $D$ must be of the form
\be
e^\lambda_{\bldr\blds}
\mapsto
d(\bldr)e^\lambda_{\bldr\blds},\quad \lambda\in\Lambda, \bldr,\blds\in\Gamma_\lambda.
\ee

From (\ref{eq:left_mult}) we now have
\bea
\lefteqn{[D,\pi(u_{ij})]e^\lambda_{\bldr\blds}} \nonumber \\
&=&
\sum_{\genfrac{}{}{0pt}{}{M\in\bbm_i, N\in\bbm_j}{m_1=n_1}}\sgn(M)\sgn(N) (d(M(\bldr))-d(\bldr))
    q^{C(i,\bldr,M)+C(j,\blds,N)+A(M)+K(M)+B(N)} 
    e^\mu_{M(\bldr)N(\blds)}.\nonumber \\
  &&  \label{eq:bdd_comm}
\eea
Therefore the condition for boundedness of commutators reads
as follows:
\be \label{eq:eqbdd1}
|(d(M(\bldr))-d(\bldr))q^{C(i,\bldr,M)+C(j,\blds,N)} |<c,
\ee
where $c$ is independent of $i$, $j$, $\lambda$, $\mu$, $\bldr$, $\blds$, $M$ and $N$.
Choosing $j$, $\blds$ and $N$ suitably, one can ensure that
(\ref{eq:eqbdd1}) implies the following:
\be\label{eq:eqbdd4}
|(d(M(\bldr))-d(\bldr)) |<c q^{-C(i,\bldr,M)}.
\ee
It follows from~(\ref{eq:bdd_comm}) that this condition is also sufficient for
the boundedness of the commutators $[D, u_{ij}]$.

Let us next form a graph $\mathcal{G}_c$  by taking the vertex set to be
\[
\Gamma:= \mbox{ the set of all GT tableaux }= \cup_\lambda \Gamma_\lambda,
\]
and connecting two elements $\bldr$ and $\bldr'$ if
$|d(\bldr)-d(\bldr')|<c$.
We will assume the existence of a partition
$(\Gamma^+,\Gamma^-)$
that has the following important property:
there does not exist infinite number of disjoint
paths each going from a point in $\Gamma^+$
to a point in $\Gamma^-$ (Here `disjoint paths' mean paths
for which the set of vertices of one does not
intersect the set of vertices of the other).
We will describe this by simply saying that the partition
$(\Gamma^+,\Gamma^-)$ does not admit any \textbf{infinite ladder}.
Existence of such a partition can be seen by looking at the sets
$\{\bldr: d(\bldr)>0\}$ and  $\{\bldr: d(\bldr)<0\}$ and exploiting the fact
that $D$ has compact resolvent.
For any subset $F$ of $\Gamma$, we will denote by $F^\pm$ the sets
$F\cap \Gamma^\pm$.
Our next job is to study this graph in more detail
using the boundedness conditions above.
Let us start with  a few definitions and notations.
By an  \textbf{elementary move}, we will mean a map $M$
from some subset of
$\Gamma$ to $\Gamma$ such that $\gamma$ and $M(\gamma)$
are connected by an edge.
A \textbf{move} will mean a composition of a finite number of
elementary moves.
If $M_1$ and $M_2$ are two moves, $M_1M_2$ and $M_2M_1$ will
in general be different.
For a family of moves $M_1, M_2,\ldots, M_r$,
we will denote by
$\overrightarrow{\prod}_{{j=1}}^{r}M_j$
the move $M_1M_2\ldots M_r$,
and by
$\overleftarrow{\prod}_{j=1}^{r}M_{j}$
the move $M_r\ldots M_2M_1$.
For a non negative integer $n$ and a move $M$, we will denote
by $M^n$ the move obtained by applying $M$ successively $n$ times.
The following families of moves will
be particularly useful to us:
\[
M_{ik}=(i,i-1,\ldots,i-k+1)\in\bbm_k,\quad
N_{ik}=(\underbrace{i+1,\ldots,i+1}_{\mbox{$k$}},
       i,i,\ldots,i)\in\bbm_{\ell+2-i}.
\]
For describing a path in our graph, we will
often use phrases like `apply the move
$\overrightarrow{\prod}_{{j=1}}^{k}M_j$
to go from $\bldr$ to $\blds$'. This will refer
to the path given by
\[
\Bigl(\bldr,\, M_k(\bldr),
   M_{k-1}M_k(\bldr),\,\ldots,\,M_1M_2\ldots
    M_k(\bldr)=\blds\Bigr).
\]
The following lemma will be very useful in the next
two sections.
\blmma\label{freemove}
Let $N_{jk}$ and $M_{ik}$ be the moves defined above. Then
\begin{enumerate}
\item  $|d(\bldr)-d(N_{j0}(\bldr))|\leq c$,
\item    $|d(\bldr)-d(M_{ik}(\bldr))|\leq
  cq^{-\sum_{a=1}^{k-1}H_{a,i+1-a}-\sum_{b=i}^{\ell}H_{k,b+k-1}}$.  In
  particular, if $H_{a,i+1-a}(\bldr)=0$ for $1\leq a\leq k-1$ and
  $H_{k,b+k-1}(\bldr)=0$ for $i\leq b\leq \ell$, then
  $|d(\bldr)-d(M_{ik}(\bldr))|\leq c$.
\end{enumerate}
\elmma
\prf
Direct consequence of~(\ref{eq:eqbdd4}).
\qed


Next, we will derive a precise estimate of the singular values of $D$.
The main ingredients in the proof are
the finiteness of exactly one of the sets $F^+$ and $F^-$
for appropriately chosen subsets $F$ of $\Gamma$.
General form of the argument for proving this will be
as follows:
for a carefully chosen coordinate $C$ (in the present case, $C$
would be one of the $V_{a1}$'s or $H_{ab}$'s), a sweepout argument
will show that any $\gamma$ can be connected by a path,
throughout which $C(\cdot)$ remains constant, to another point $\gamma'$
for which $C(\gamma')=C(\gamma)$ and all other coordinates of $\gamma'$
are zero.
This would help connect any two points $\gamma$ and $\delta$ by a path
such that $C(\cdot)$ would lie between $C(\gamma)$ and $C(\delta)$
on the path. This would finally result in the finiteness
of at least one (and hence exactly one) of $C(F^+)$ and $C(F^-)$.
Next, assuming one of these, say $C(F^-)$ is finite,
one shows that for any other coordinate $C'$,
$C'(F^-)$ is also finite.
This is done as follows. If $C'(F^-)$ is infinite, one chooses
elements $y_n\in F^-$ with
$C'(y_n)<C'(y_{n+1})$ for all $n$.
Now starting at each $y_n$, produce paths
keeping the $C'$-coordinate constant and taking the
$C$-coordinate above the plane $C(\cdot)=K$, where $C(F^-)\seq [-K,K]$.
This will produce an infinite ladder.
The argument is explained in the following diagram.\\[3ex]
\hspace*{60pt}
\setlength{\unitlength}{0.00041667in}
\begingroup\makeatletter\ifx\SetFigFont\undefined%
\gdef\SetFigFont#1#2#3#4#5{%
  \reset@font\fontsize{#1}{#2pt}%
  \fontfamily{#3}\fontseries{#4}\fontshape{#5}%
  \selectfont}%
\fi\endgroup%
{\renewcommand{\dashlinestretch}{30}
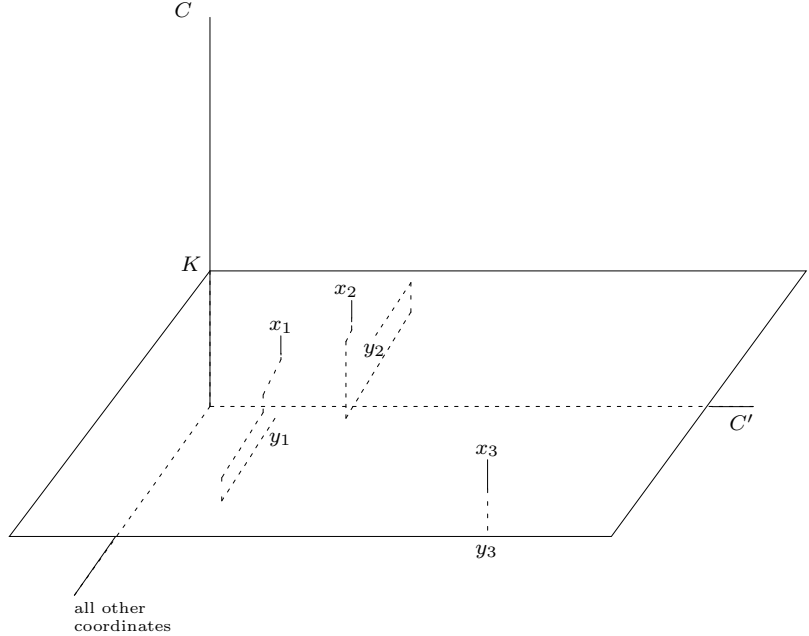
\begin{picture}(10149,7971)(0,-10)
\path(2562,2850)(2562,7800)
\dashline{60.000}(2562,2850)(837,450)
\dashline{60.000}(9462,2850)(2562,2850)(2562,4575)
\dashline{60.000}(4287,3675)(4287,2700)(5112,4050)
        (5112,4425)(4662,3750)
\path(1362,1200)(837,450)
\path(8937,2850)(9462,2850)
\path(2562,4575)(12,1200)(7662,1200)
        (10137,4575)(2562,4575)
\path(3462,3750)(3462,3525)
\dashline{60.000}(3462,3525)(3462,3450)(3237,3000)
        (3237,2775)(2712,1950)(2712,1650)(3387,2700)
\dashline{60.000}(4287,3675)(4362,3825)(4362,3975)
\path(4362,3975)(4362,4200)
\dashline{60.000}(6087,1275)(6087,1800)
\path(6087,1800)(6087,2175)
\put(2112,7800){\makebox(0,0)[lb]{\smash{{{\SetFigFont{8}{9.6}{\rmdefault}{\mddefault}{\updefault}$C$}}}}}
\put(2187,4575){\makebox(0,0)[lb]{\smash{{{\SetFigFont{8}{9.6}{\rmdefault}{\mddefault}{\updefault}$K$}}}}}
\put(837,225){\makebox(0,0)[lb]{\smash{{{\SetFigFont{6}{7.2}{\rmdefault}{\mddefault}{\updefault}all other}}}}}
\put(837,0){\makebox(0,0)[lb]{\smash{{{\SetFigFont{6}{7.2}{\rmdefault}{\mddefault}{\updefault}coordinates}}}}}
\put(9162,2550){\makebox(0,0)[lb]{\smash{{{\SetFigFont{8}{9.6}{\rmdefault}{\mddefault}{\updefault}$C'$}}}}}
\put(4512,3525){\makebox(0,0)[lb]{\smash{{{\SetFigFont{8}{9.6}{\rmdefault}{\mddefault}{\updefault}$y_2$}}}}}
\put(3312,2400){\makebox(0,0)[lb]{\smash{{{\SetFigFont{8}{9.6}{\rmdefault}{\mddefault}{\updefault}$y_1$}}}}}
\put(3312,3825){\makebox(0,0)[lb]{\smash{{{\SetFigFont{8}{9.6}{\rmdefault}{\mddefault}{\updefault}$x_1$}}}}}
\put(5937,975){\makebox(0,0)[lb]{\smash{{{\SetFigFont{8}{9.6}{\rmdefault}{\mddefault}{\updefault}$y_3$}}}}}
\put(5937,2250){\makebox(0,0)[lb]{\smash{{{\SetFigFont{8}{9.6}{\rmdefault}{\mddefault}{\updefault}$x_3$}}}}}
\put(4137,4275){\makebox(0,0)[lb]{\smash{{{\SetFigFont{8}{9.6}{\rmdefault}{\mddefault}{\updefault}$x_2$}}}}}
\end{picture}
}\\[2ex]

Our next job is to define an important class of subsets of $\Gamma$.
Observe that lemma~\ref{freemove}
tells us that for any $\bldr$ and any $j$, the points
$\bldr$ and $N_{j0}(\bldr)$ are connected by an edge,
whenever $N_{j0}(\bldr)$ is a GT tableaux.
Let $\bldr$ be an element of $\Gamma$. Define the
\textbf{free plane passing through $\bldr$} to be  the minimal
subset of $\Gamma$ that contains $\bldr$
and is closed under application of the moves $N_{j0}$.
We will denote this set by $\scrf_\bldr$.
The following is an easy consequence of this definition.
\blmma \label{freecriterion}
Let $\bldr$ and $\blds$ be two GT tableaux. Then
$\blds\in \scrf_\bldr$ if and only if
$V_{a,1}(\bldr)=V_{a,1}(\blds)$ for all $a$ and for each $b$, the difference
$H_{a,b}(\bldr)-H_{a,b}(\blds)$ is independent of $a$.
\elmma
\bcrlre \label{freedisjt}
Let $\bldr,\blds\in\Gamma$. Then either $\scrf_\bldr=\scrf_\blds$ or
$\scrf_\bldr\cap \scrf_\blds=\phi$.
\ecrlre

Let $\bldr\in\Gamma$.
For $1\leq j\leq \ell+1$, define
$a_j$ to be an integer such that $H_{a_j,j}(\bldr)=\min_i H_{ij}(\bldr)$.
Note three things here:\\
1. definition of $a_j$ depends on $\bldr$,\\
2. for a given $j$ and given $\bldr$, $a_j$ need not be unique, and\\
3. if $\blds\in\scrf_\bldr$, then for each $j$, the set of
$k$'s for which $H_{kj}(\blds)=\min_i H_{ij}(\blds)$ is same
as the set of all $k$'s for which $H_{kj}(\bldr)=\min_i H_{ij}(\bldr)$.
Therefore, the $a_j$'s can be chosen in a manner such that
they remain the same for all elements lying on a given free plane.

\blmma\label{sweep1}
Let $\blds\in \scrf_\bldr$. Let $\blds'$ be another GT tableaux
given by
\[
V_{a1}(\blds')=V_{a1}(\blds) \mbox{ and }
H_{a1}(\blds')=H_{a1}(\blds) \mbox{ for all }a,\quad
 H_{a_b,b}(\blds')=0 \mbox{ for all }b>1,
\]
where the $a_j$'s are as defined above.
Then there is a path in $\scrf_\bldr$ from $\blds$ to $\blds'$
such that $H_{11}(\cdot)$ remains constant throughout this path.
\elmma
\prf
Let $c_b:=\sum_{j=2}^{\ell+2-b} H_{a_j,j}(\blds)$. Apply the move
$\overrightarrow{\prod}_{{b=2}}^{\ell}
   N_{\ell+3-b,0}^{c_b}$.\qed

The following diagram will help explain the steps involved
in the above proof in the case where $\bldr$ is the constant
tableaux.\\[2ex]
\def\labelstyle{\scriptstyle}
\xymatrix@C=.6pt@R=.6pt{
    \cdot\ar@{}[r] & & \cdot\ar@{}[r] && \odot\ar@{.}[d] && \cdot&& \cdot&\\
    0 & a &  & b & & c && d &\\
    \cdot\ar@{}[r] & & \cdot\ar@{}[r] && \cdot\ar@{.}[u]\ar@{.}[d] && \cdot&\\
    0 & a &  & b & & c &\\
    \cdot\ar@{}[r] & & \cdot\ar@{}[r] && \odot\ar@{.}[u]&\\
    0 & a && b &  \\
    \cdot\ar@{}[r] && \cdot&\\
    0 &a & \\
    \cdot&}
\hspace{-2em} \xymatrix@C=20pt@R=12pt{&\\&\\ \ar@{->}[r]^{N_{30}^b}&\\}\hspace{.3em}
\xymatrix@C=.6pt@R=.6pt{
    \cdot\ar@{}[r] & & \cdot\ar@{}[r] && \cdot && \odot\ar@{.}[d]&& \cdot&\\
    0 & a &  & 0 & & b+c && d &\\
    \cdot\ar@{}[r] & & \cdot\ar@{}[r] && \cdot && \odot\ar@{.}[u]&\\
    0 & a &  & 0 & & b+c &\\
    \cdot\ar@{}[r] & & \cdot\ar@{}[r] && \cdot&\\
    0 & a && 0 &  \\
    \cdot\ar@{}[r] && \cdot&\\
    0 &a & \\
    \cdot&}
\hspace{-2em} \xymatrix@C=20pt@R=12pt{&\\&\\ \ar@{->}[r]^{N_{40}^{b+c}}&\\}\hspace{.3em}
\xymatrix@C=.6pt@R=.6pt{
    \cdot\ar@{}[r] & & \cdot\ar@{}[r] && \cdot && \cdot& \odot&\\
    0 & a &  & 0 & & 0 && b+c+d &\\
    \cdot\ar@{}[r] & & \cdot\ar@{}[r] && \cdot && \cdot&\\
    0 & a &  & 0 & & 0 &\\
    \cdot\ar@{}[r] & & \cdot\ar@{}[r] && \cdot&\\
    0 & a && 0 &  \\
    \cdot\ar@{}[r] && \cdot&\\
    0 &a & \\
    \cdot&}\\[2ex]
  \hspace*{12em}  \xymatrix@C=20pt@R=12pt{&\\&\\ \ar@{->}[r]^{N_{50}^{b+c+d}}&\\}\hspace{.3em}
\xymatrix@C=.6pt@R=.6pt{
    \cdot\ar@{}[r] & & \cdot\ar@{}[r] && \cdot && \cdot&& \cdot&\\
    0 & a &  & 0 & & 0 && 0 &\\
    \cdot\ar@{}[r] & & \cdot\ar@{}[r] && \cdot && \cdot&\\
    0 & a &  & 0 & & 0 &\\
    \cdot\ar@{}[r] & & \cdot\ar@{}[r] && \cdot&\\
    0 & a && 0 &  \\
    \cdot\ar@{}[r] && \cdot&\\
    0 &a & \\
    \cdot&}\\
A dotted line joining two circled dots signifies a move that
increases the $r_{ij}$'s lying on the dotted line by one.
Where there is one circled dot and no dotted line, it means
one applies the move that raises the $r_{ij}$ corresponding to
the circled dot by one.

\bppsn\label{signfree1}
Let $\bldr$ be a GT tableaux.
Then either $\scrf_{\bldr}^+$ is finite or $\scrf_{\bldr}^-$ is finite.
\eppsn
\prf
Suppose, if possible, both $H_{11}(\scrf_{\bldr}^+)$ and $H_{11}(\scrf_{\bldr}^-)$
are infinite. Then there exist two sequences of elements
$\bldr_n$ and $\blds_n$ with $\bldr_n\in \scrf_\bldr^+$
and $\blds_n\in \scrf_\bldr^-$,
such that
\[
H_{11}(\bldr_1)<H_{11}(\blds_1)<H_{11}(\bldr_2)<H_{11}(\blds_2)<\cdots.
\]
Now starting from $\bldr_n$,  employ the forgoing lemma
to reach a point $\bldr'_n\in\scrf_{\bldr}$ for which
\[
V_{a1}(\bldr'_n)=V_{a1}(\bldr_n) \mbox{ and }
H_{a1}(\bldr'_n)=H_{a1}(\bldr_n) \mbox{ for all }a,\quad
 H_{a_b,b}(\bldr'_n)=0 \mbox{ for all }b>1.
\]
Similarly, start at $\blds_n$ and go to
a point $\blds'_n\in\scrf_{\bldr}$ for which
\[
V_{a1}(\blds'_n)=V_{a1}(\blds_n) \mbox{ and }
H_{a1}(\blds'_n)=H_{a1}(\blds_n) \mbox{ for all }a,\quad
 H_{a_b,b}(\blds'_n)=0 \mbox{ for all }b>1.
\]
Now use the move $N_{10}$ to get to $\blds'_n$ from $\bldr'_n$.
The paths thus constructed are all disjoint, because
for the path from $\bldr_n$ to $\blds_n$, the
$H_{11}$ coordinate lies between
$H_{11}(\bldr_n)$ and $H_{11}(\blds_n)$.
This means $(\scrf_\bldr^+, \scrf_\bldr^-)$ admits an infinite ladder,
which can not happen.
So one of the sets  $H_{11}(\scrf_\bldr^+)$ and $H_{11}(\scrf_\bldr^-)$
must be finite. Let us assume that $H_{11}(\scrf_\bldr^-)$ is finite.

Let us next show that for any $b>1$, $H_{ab}(\scrf_\bldr^-)$ is finite.
Let $K$ be an integer such that $H_{11}(\blds)<K$ for all $\blds\in \scrf_\bldr^-$.
If $H_{ab}(\scrf_\bldr^-)$ was infinite, there would exist elements
$\bldr_n\in \scrf_\bldr^-$ such that
\[
H_{ab}(\bldr_1)<H_{ab}(\bldr_2)<\cdots.
\]
Now start at $\bldr_n$ and employ the move $N_{10}$ successively
$K$ times to reach a point in
$\scrf_\bldr^+=\scrf_\bldr\backslash\scrf_\bldr^-$.
These paths will all be disjoint, as throughout the path,
$H_{ab}$ remains fixed.

Since the coordinates $(H_{11},H_{12},\ldots,H_{1,\ell})$
completely specify a point in $\scrf_\bldr$, it follows that
$\scrf_\bldr^-$ is finite.
\qed

Next we need a set that can be used for a proper
indexing of the free planes.
Such a set will be called a complementary axis.

\bdfn\rm
A subset $\scrc $ of $\Gamma$ is called a \textbf{complementary axis} if
\begin{enumerate}
\item $\cup_{\bldr\in \scrc }\scrf_\bldr =\Gamma$,
\item if $\bldr,\blds\in \scrc $, and $\bldr\neq \blds$, then
   $\scrf_\bldr$ and $\scrf_\blds$ are disjoint.
\end{enumerate}
\edfn

Let us next give a choice of a complementary axis.
\bthm \label{compl}
Define
\[
\scrc =\{\bldr\in \Gamma: \Pi_{a=1}^{\ell+1-b} H_{ab}(\bldr)=0
           \mbox{ for } 1\leq b\leq \ell\}.
\]
The set $\scrc $ defined above is a complementary axis.
\ethm
\prf
Let $\blds\in\Gamma$.
A sweepout argument almost identical to that used in
lemma~\ref{sweep1} (application of the move
$\overrightarrow{\prod}_{{b=1}}^\ell
   N_{\ell+2-b,0}^{\sum_{j=1}^{\ell+1-b} H_{a_j,j}(\blds)}$ )
will connect $\blds$ to another element $\blds'$
for which $H_{a_b,b}(\blds')=0$ for $1\leq b\leq\ell$
by a path that lies entirely on $\scrf_\blds$.
Clearly, $\blds'\in\scrc$. Since $\blds'\in\scrf_\blds$,
by corollary~\ref{freedisjt}, $\blds\in\scrf_{\blds'}$.

It remains to show that if $\bldr$ and $\blds$ are two distinct elements of
$\scrc$, then $\blds\not\in\scrf_\bldr$.
Since $\bldr\neq\blds$, there exist two integers  $a$ and $b$,
$1\leq b\leq \ell$ and $1\leq a\leq \ell+2-b$, such that
$H_{ab}(\bldr)\neq H_{ab}(\blds)$.
Observe that $H_{1\ell}(\cdot)$ must be zero for both,
as they are members of $\scrc$. So $b$ can not be $\ell$ here.
Next we will produce two integers $i$ and $j$ such that
the differences $H_{ib}(\bldr)-H_{ib}(\blds)$
and $H_{jb}(\bldr)-H_{jb}(\blds)$ are distinct.
If there is an integer $k$ for which
$H_{kb}(\bldr)=H_{kb}(\blds)=0$, then take $i=a$, $j=k$.
If not, there would exist two integers $i$ and $j$ such that
$H_{ib}(\bldr)=0$, $H_{ib}(\blds)>0$
and
$H_{jb}(\bldr)>0$, $H_{jb}(\blds)=0$.
Take these $i$ and $j$.
Since $H_{ib}(\bldr)-H_{ib}(\blds)$
and $H_{jb}(\bldr)-H_{jb}(\blds)$ are distinct,
by lemma~\ref{freecriterion}, $\bldr$ and $\blds$ can not lie
on the same free plane.
\qed

\blmma \label{sweep2}
Let $\bldr$ be a GT tableaux. Let $\blds$ be the GT tableaux
defined by the prescription
\[
V_{a1}(\blds)=V_{a1}(\bldr)\mbox{ for all }a,\quad
 H_{ab}(\blds)=H_{ab}(\bldr) \mbox{ for all }a\geq 2,\mbox{ for all }b,\quad
 H_{1,b}(\blds)=0 \mbox{ for all }b.
\]
Then there is a path from $\bldr$ to $\blds$ such that
$V_{a1}(\cdot)$ remains constant throughout the path.
\elmma
\prf
Apply  the move
$\overrightarrow{\prod}_{b=1}^\ell M_{b+1,1}^{H_{1,b}(\bldr)}$.\qed

The above lemma is actually the first step in the following
slightly more general sweepout algorithm.

\blmma \label{sweep3}
Let $\bldr$ be a GT tableaux. Let $\blds$ be the GT tableaux
defined by the prescription
\[
V_{11}(\blds)=V_{11}(\bldr),\quad
V_{a1}(\blds)=0\mbox{ for all }a>1,\quad
 H_{ab}(\blds)=0 \mbox{ for all }a,b.
\]
Then there is a path from $\bldr$ to $\blds$ such that
$V_{11}(\cdot)$ remains constant throughout the path.
\elmma
\prf
Apply successively the moves
\[
\overrightarrow{\prod}_{b=1}^\ell M_{b+1,1}^{H_{1,b}(\bldr)},\quad
\overrightarrow{\prod}_{{b=1}}^{\ell-1} M_{b+2,2}^{H_{2,b}(\bldr)},\quad
\ldots,\quad
 M_{\ell+1,\ell}^{H_{\ell,1}(\bldr)},
\]
followed by
\be\label{movseq}
M_{33}^{V_{21}(\bldr)},\quad  M_{44}^{V_{21}(\bldr)+V_{31}(\bldr)},\quad
\ldots,\quad  M_{\ell+1,\ell+1}^{\sum_{a=2}^\ell V_{a1}(\bldr)}.
\ee
\qed

The following diagram
will help explain the procedure described above in a simple case.\\
\def\labelstyle{\scriptstyle}
\xymatrix@C=10pt@R=8pt{
    \cdot\ar@{}[r]\ar@{}[d]|\star      & \cdot\ar@{}[r] & \cdot\ar@{}[r] & \odot&\\
    \cdot\ar@{}[d]|\star\ar@{}[ur]|\star& \cdot\ar@{}[ur]|\star & \cdot\ar@{}[ur]|p &  \\
    \cdot\ar@{}[d]|\star\ar@{}[ur]|\star& \cdot\ar@{}[ur]|\star & \\
     \cdot\ar@{}[ur]|\star& }
 \hspace{-1em} \xymatrix@C=20pt@R=12pt{&\\ \ar@{->}[r]^{M_{41}^p}&\\}\hspace{.5em}
\xymatrix@C=10pt@R=8pt{
    \cdot\ar@{}[r]\ar@{}[d]|\star      & \cdot\ar@{}[r] & \odot\ar@{}[r] & \cdot&\\
    \cdot\ar@{}[d]|\star\ar@{}[ur]|\star& \cdot\ar@{}[ur]|q & \cdot\ar@{}[ur]|0 &  \\
    \cdot\ar@{}[d]|\star\ar@{}[ur]|\star& \cdot\ar@{}[ur]|\star & \\
     \cdot\ar@{}[ur]|\star&  }
  \hspace{-1em}  \xymatrix@C=20pt@R=12pt{&\\ \ar@{->}[r]^{M_{31}^q}&\\}\hspace{.5em}
\xymatrix@C=10pt@R=8pt{
    \cdot\ar@{}[r]\ar@{}[d]|\star      & \odot\ar@{}[r] & \cdot\ar@{}[r] & \cdot&\\
    \cdot\ar@{}[d]|\star\ar@{}[ur]|r& \cdot\ar@{}[ur]|0 & \cdot\ar@{}[ur]|0 &  \\
    \cdot\ar@{}[d]|\star\ar@{}[ur]|\star& \cdot\ar@{}[ur]|\star & \\
     \cdot\ar@{}[ur]|\star& \\
     }
  \hspace{-1em}  \xymatrix@C=20pt@R=12pt{&\\ \ar@{->}[r]^{M_{21}^r}&\\}\hspace{.5em}
\xymatrix@C=10pt@R=8pt{
    \cdot\ar@{}[r]\ar@{}[d]|\star      & \cdot\ar@{}[r] & \cdot\ar@{}[r] & \odot&\\
    \cdot\ar@{}[d]|\star\ar@{}[ur]|0& \cdot\ar@{}[ur]|0 & \odot\ar@{.}[ur]|0 &  \\
    \cdot\ar@{}[d]|\star\ar@{}[ur]|\star& \cdot\ar@{}[ur]|s & \\
     \cdot\ar@{}[ur]|\star& \\
     }\\[2ex]
\hspace*{60pt}
  \xymatrix@C=20pt@R=12pt{&\\ \ar@{->}[r]^{M_{42}^s}&\\}\hspace{.5em}
\xymatrix@C=10pt@R=8pt{
    \cdot\ar@{}[r]\ar@{}[d]|\star      & \cdot\ar@{}[r] & \odot\ar@{}[r] & \cdot&\\
    \cdot\ar@{}[d]|\star\ar@{}[ur]|0& \odot\ar@{.}[ur]|0 & \cdot\ar@{}[ur]|0 &  \\
    \cdot\ar@{}[d]|\star\ar@{}[ur]|t& \cdot\ar@{}[ur]|0 & \\
     \cdot\ar@{}[ur]|\star& \\
     }
  \hspace{-1em}  \xymatrix@C=20pt@R=12pt{&\\ \ar@{->}[r]^{M_{32}^t}&\\}\hspace{.5em}
\xymatrix@C=10pt@R=8pt{
    \cdot\ar@{}[r]\ar@{}[d]|\star      & \cdot\ar@{}[r] & \cdot\ar@{}[r] & \odot&\\
    \cdot\ar@{}[d]|\star\ar@{}[ur]|0& \cdot\ar@{}[ur]|0 & \cdot\ar@{.}[ur]|0 &&  \\
    \cdot\ar@{}[d]|\star\ar@{}[ur]|0& \odot\ar@{.}[ur]|0 &&& \\
     \cdot\ar@{}[ur]|u& &&&\\
     }
  \hspace{-1em}  \xymatrix@C=20pt@R=12pt{&\\ \ar@{->}[r]^{M_{43}^u}&\\}\hspace{.5em}
\xymatrix@C=10pt@R=8pt{
    \cdot\ar@{}[r]\ar@{}[d]|\star      & \cdot\ar@{}[r] & \odot\ar@{}[r] & \cdot&\\
    \cdot\ar@{}[d]|v\ar@{}[ur]|0& \cdot\ar@{.}[ur]|0 & \cdot\ar@{}[ur]|0 &  \\
    \odot\ar@{}[d]|\star\ar@{.}[ur]|0& \cdot\ar@{}[ur]|0 & \\
     \cdot\ar@{}[ur]|0& \\
     }\\[2ex]
\hspace*{60pt}
  \xymatrix@C=20pt@R=12pt{&\\ \ar@{->}[r]^{M_{33}^v}&\\}\hspace{.5em}
\xymatrix@C=10pt@R=8pt{
    \cdot\ar@{}[r]\ar@{}[d]|\star      & \cdot\ar@{}[r] & \cdot\ar@{}[r] & \odot&\\
    \cdot\ar@{}[d]|0\ar@{}[ur]|0& \cdot\ar@{}[ur]|0 & \cdot\ar@{.}[ur]|0 &  \\
    \cdot\ar@{}[d]|w\ar@{}[ur]|0& \cdot\ar@{.}[ur]|0 & \\
     \odot\ar@{.}[ur]|0& \\
     }
  \xymatrix@C=20pt@R=12pt{&\\ \ar@{->}[r]^{M_{44}^w}&\\}\hspace{.5em}
\xymatrix@C=10pt@R=8pt{
    \cdot\ar@{}[r]\ar@{}[d]|\star      & \cdot\ar@{}[r] & \cdot\ar@{}[r] & \cdot&\\
    \cdot\ar@{}[d]|0\ar@{}[ur]|0& \cdot\ar@{}[ur]|0 & \cdot\ar@{}[ur]|0 &  \\
    \cdot\ar@{}[d]|0\ar@{}[ur]|0& \cdot\ar@{}[ur]|0 & \\
     \cdot\ar@{}[ur]|0& \\
     }

\bcrlre\label{growth5}
$|d(\bldr)|=O(r_{11})$.
\ecrlre
\prf
If one employs the sequence of moves
\[
M_{22}^{V_{11}(\bldr)},\quad M_{33}^{V_{11}(\bldr)+V_{21}(\bldr)},\quad
\ldots,\quad  M_{\ell+1,\ell+1}^{\sum_{a=1}^\ell V_{a1}(\bldr)}
\]
instead of the sequence given in (\ref{movseq}),
one would reach the constant (or zero) tableaux.
Total length of this path from $\bldr$ to the zero tableaux
is
\[
\sum_{a=1}^\ell\sum_{b=1}^{\ell+1-a}H_{ab}(\bldr)
   + \sum_{b=1}^\ell\sum_{a=1}^b V_{a1}(\bldr),
\]
which can easily be shown to be bounded by $\ell r_{11}$.
\qed


\bthm\label{th:summability}
Let $\widetilde{D}$  be the following operator:
\be
\widetilde{D}: e^\lambda_{\bldr,\blds}\mapsto r_{11}e^\lambda_{\bldr,\blds}
\ee
Then $(\cla,\clh,\widetilde{D})$ is an equivariant  spectral triple.

Moreover, $\widetilde{D}$ is not $p$-summable if $p\leq\ell(\ell+2)$,
but is $p$-summable for all $p>\ell(\ell+2)$.
\ethm
\prf
Boundedness of commutators with algebra elements
follow from the observation that
$|d(\bldr)-d(M(\bldr)|\leq 1$ and hence
equation~(\ref{eq:eqbdd4}) is satisfied.

Let $\delta_n$ denote the dimension of the eigenspace
of $\widetilde{D}$ corresponding to the eigenvalue $n$.

Observe that the number of Young tableux
$\lambda=(\lambda_1,\ldots,\lambda_\ell,\lambda_{\ell+1})$
with
$n=\lambda_1\geq \lambda_2\geq \ldots \lambda_\ell\geq \lambda_{\ell+1}=0$
 is
\[
\sum_{i_1=0}^{n}\sum_{i_2=0}^{i_1}\ldots\sum_{i_{\ell-1}=0}^{i_{\ell-2}}1
=
\mbox{polynomial in $n$ of degree $\ell-1$}.
\]
Thus the number of such Young tableaux is $O(n^{\ell-1})$.

Next, let
$\lambda:n=\lambda_1\geq \lambda_2\geq\ldots\geq \lambda_{\ell}\geq 0$
be an Young tableaux, and let $V_\lambda$ be the space
carrying the irreducible representation parametrized by
$\lambda$. Then by Weyl dimension formula,
\bean
\mbox{dim}\, V_\lambda &=&
 \prod_{1\leq i<j\leq \ell+1}
  \frac{(\lambda_{i}-\lambda_{i+1})+
      \ldots + (\lambda_{j-1}-\lambda_{j})+j-i}{j-i}\\
  &=& \prod_{1\leq i<j\leq \ell+1}
  \frac{\lambda_{i}-\lambda_{j}+j-i}{j-i}\\
  &\leq& (n+1)^{\frac{\ell(\ell+1)}{2}}.
\eean
Thus the dimension of an irreducible representation corresponding to
a Young tableaux
\[
n=\lambda_1\geq \lambda_2\geq \ldots \lambda_\ell\geq \lambda_{\ell+1}=0
\]
is $O(n^{\half\ell(\ell+1)})$.

Using the two observations above, 
it follows that 
\[
\delta_n\leq C n^{\ell-1}\left(n^{\half\ell(\ell+1)}\right)^2=C n^{\ell(\ell+2)-1},
\]
where $C$ denotes a generic constant. 


This implies that for $p>\ell(\ell+2)$, one has
\[
\mbox{Trace\,}|\widetilde{D}|^{-p}\leq 
\sum_n n^{-p}n^{\ell(\ell+2)-1}=\sum_n\frac{1}{n^{1+p-\ell(\ell+2)}}<\infty,
\]
i.e.\ $\widetilde{D}$ is $p$-summable.

Next, let us take an $\epsilon\in (0,\frac{1}{4\ell})$.
Then for large enough $n$, the number of Young tableux
$\lambda=(\lambda_1,\ldots,\lambda_\ell,\lambda_{\ell+1})$
with
\[
\lambda_1=n,\quad \left|\lambda_2-\left(1-\frac{1}{\ell}\right)n\right|<\epsilon n, \quad
  \left|\lambda_3-\left(1-\frac{2}{\ell}\right)n\right|<\epsilon n, \quad \ldots \quad
  \left|\lambda_\ell-\frac{1}{\ell}n\right|<\epsilon n, \quad \lambda_{\ell+1}=0
\]
is of the order $n^{\ell-1}$.
For each such $\lambda$ and for $1\leq i < j \leq \ell+1$,
one has
$\lambda_i-\lambda_j> \left(\frac{j-i}{\ell}-2\epsilon\right)n$, so that
\[
\frac{\lambda_i-\lambda_j}{j-i}>\left(\frac{1}{\ell}-\frac{2\epsilon}{j-i}\right)n > \frac{n}{2\ell}.
\]
Therefore
\[
\mbox{dim}\, V_\lambda =
 \prod_{1\leq i<j\leq \ell+1}
  \frac{\lambda_{i}-\lambda_{j}+j-i}{j-i}
      > \left(\frac{n}{2\ell}+1\right)^{\half\ell(\ell+1)} > C n^{\half\ell(\ell+1)},
\]
$C$ being a generic constant.
It now follows that
\[
\delta_n\geq C n^{\ell-1}\left(n^{\half\ell(\ell+1)}\right)^2=C n^{\ell(\ell+2)-1}.
\]
Hence for $p\leq\ell(\ell+2)$, we have
\[
\mbox{Trace\,}|\widetilde{D}|^{-p}\geq \sum_n n^{-p}n^{\ell(\ell+2)-1}
  =\sum_n \frac{1}{n^{1+p-\ell(\ell+2)}}=\infty,
\]
i.e.\ $\widetilde{D}$ can not be $p$-summable.
\qed

As an important consequence of  corollary~\ref{growth5} and the above theorem, 
we now derive the following.

\bthm
The spectral dimension of the quantum group $SU_q(\ell+1)$ is $\ell(\ell+2)$.
\ethm
\prf
If $(G, L_2(G), D)$ is an equivariant
spectral triple, then by corollary~\ref{growth5}, the singular values of $D$ grow slower
than those of (a scalar multiple of) $\widetilde{D}$. Since $\widetilde{D}$ is not $p$-summable
for $p\leq \ell(\ell+2)$, the operator $D$ also can not be $p$-summable
for $p\leq \ell(\ell+2)$. On the other hand, from theorem~\ref{th:summability}
we know that $(G, L_2(G),\widetilde{D})$ is an equivariant spectral triple that
is $p$-summable for all $p>\ell(\ell+2)$.
Therefore the result follows.
\qed

\section{$SU_q(\ell+1)$ action on $S_q^{2\ell+1}$}
The $C^*$-algebra $A_\ell\equiv C(S_q^{2\ell+1})$ of the quantum
sphere $S_q^{2\ell+1}$
is the universal $C^*$-algebra generated by
elements
$z_1, z_2,\ldots, z_{\ell+1}$
satisfying the following relations (see~\cite{hon-szy-2002a}):
\bean
z_i z_j & =& qz_j z_i,\qquad 1\leq j<i\leq \ell+1,\\
z_i^* z_j & =& q z_j z_i^* ,\qquad 1\leq i\neq j\leq \ell+1,\\
z_i z_i^* - z_i^* z_i +
(1-q^{2})\sum_{k>i} z_k z_k^* &=& 0,\qquad \hspace{2em}1\leq i\leq \ell+1,\\
\sum_{i=1}^{\ell+1} z_i z_i^* &=& 1.
\eean

Let $u_{ij}$ denote the generating elements of the $C^*$-algebra
$C(SU_q(\ell+1))$ as in the previous section.
The map
\[
\tau(z_i)=\sum_k z_k\otimes u_{ki}^*
\]
extends to a *-homomorphism $\tau$ from $A_\ell$
into $A_\ell\otimes C(SU_q(\ell+1))$ and
obeys $(\mbox{id}\otimes\Delta)\tau = (\tau\otimes\mbox{id})\tau$.
In other words this gives an action of $SU_q(\ell+1)$ on $A_\ell$.
Equivariant spectral triples for this dynamical system were studied
in \cite{cha-pal-2008a}. As we shall see shortly, that this dynamical system
is ergodic as well as the computation of the spectral dimension of this dynamical
system is a by product of the results there.


Let us recall from  the description 
of the $L_2$ space of the sphere
sitting inside $L_2(SU_q(\ell+1))$.
Let $u^\one$ denote the fundamental unitary for $SU_q(\ell+1)$,
i.\ e.\ the irreducible unitary representation corresponding to the
Young tableaux $\one=(1,0,\ldots,0)$.
Similarly write $v^\one$ for the fundamental unitary for $SU_q(\ell)$.
Fix some bases for the corresponding representation spaces.
Then recall (\cite{cha-pal-2008a}) that $C(SU_q(\ell+1))$ is the 
$C^*$-algebra generated by the matrix
entries $\{u^\one_{ij}\}$ and $C(SU_q(\ell))$ is the
$C^*$-algebra generated by the matrix
entries $\{v^\one_{ij}\}$.
Now define $\phi$ by
\be
\phi(u^\one_{ij})=\begin{cases} I & \mbox{if $i=j=1$},\cr
        v^\one_{i-1,j-1} & \mbox{if $2\leq i,j\leq \ell+1$},\cr
        0 & \mbox{otherwise.}
       \end{cases}
\ee
Then $C(SU_q(\ell+1)\verb1\1SU_q(\ell))$ is the $C^*$-subalgebra of
$C(SU_q(\ell+1))$ generated by the entries $u_{1,j}$ for $1\leq j\leq \ell+1$.
Define $\psi:C(S_q^{2\ell+1})\rightarrow C(SU_q(\ell+1)\verb1\1SU_q(\ell))$
by 
\[
 \psi(z_i)=q^{-i+1}u^*_{1,i}.
\]
This gives an isomorphism between $C(SU_q(\ell+1)\verb1\1SU_q(\ell))$ and
$C(S_q^{2\ell+1})$,
and the following diagram commutes:
\[
\def\labelstyle{\scriptstyle}
  \xymatrix@C=23pt@R=50pt{
 C(S_q^{2\ell+1})\ar[d]_{\psi}\ar[r]_-{\tau} & C(S_q^{2\ell+1})\otimes
C(SU_q(\ell+1))\ar[d]_{\psi\otimes\id}&\\
C(SU_q(\ell+1)\backslash SU_q(\ell))\ar[r]_-{\Delta}  & C(SU_q(\ell+1)\backslash SU_q(\ell))\otimes C(SU_q(\ell+1))\\
}
\]
In other words,
$(C(S_q^{2\ell+1}), SU_q(\ell+1),\tau)$ is the quotient space $SU_q(\ell+1)\verb1\1SU_q(\ell)$.
Thus by proposition~1.9, \cite{pod-1995a}, the action we are considering
is ergodic.
Also,  note that $C(S_q^3)\cong C(SU_q(2))$ and the $SU_q(2)$-action on $S_q^3$
under this equivalence is same as the $SU_q(2)$-action on itself,
which has been covered in the previous section. Therefore we will
assume in the rest of this section that $\ell>1$.

The choice of $\psi$  makes $L_2(SU_q(\ell+1)\verb1\1SU_q(\ell))$ 
a span of certain rows of the $e_{\bldr,\blds}$'s.
To be more precise, the right regular representation $u$ of $SU_q(\ell+1)$ 
keeps the subspace $L_2(SU_q(\ell+1)\verb1\1SU_q(\ell))$ invariant, and the restriction of $u$ to
$L_2(SU_q(\ell+1)\verb1\1SU_q(\ell))$ decomposes as a direct sum of exactly one copy
of each of the irreducibles given by the young tableaux
$\lambda_{n,k}:=(n+k, k,k,\ldots, k,0)$, with $n,k\in\bbn$.
Let $\Gamma_0$ be the set of all GT tableaux $\bldr^{nk}$
given by
\[
r^{nk}_{ij}=\begin{cases} n+k & \mbox{if $i=j=1$},\cr
                0  & \mbox{if $i=1$, $j=\ell+1$},\cr
                k  & \mbox{otherwise},\end{cases}
\]
for some $n,k \in \bbn$.
Let $\Gamma_0^{nk}$ be the set of all GT tableaux with
top row $\lambda_{n,k}$.
Then the family of vectors
\[
\{e_{\bldr^{nk},\blds}: n,k\in\bbn,\, \blds\in\Gamma_0^{nk}\}
\]
form a complete
orthonormal basis for $L_2(SU_q(\ell+1)\verb1\1SU_q(\ell))$.
Thus the right regular representation $u$ restricts to the subspace
$L_2(SU_q(\ell+1)\verb1\1SU_q(\ell))$ and it also follows from the above discussion
and equation~(\ref{eq:left_mult})
that the restriction of the left multiplication to $C(SU_q(\ell+1)\verb1\1SU_q(\ell))$
keeps $L_2(SU_q(\ell+1)\verb1\1SU_q(\ell))$ invariant.
Let us denote the restriction of $u$ to $L_2(SU_q(\ell+1)\verb1\1SU_q(\ell))$
by $\hat{u}$ and the restriction of $\pi$ to $C(SU_q(\ell+1)\verb1\1SU_q(\ell))$
viewed as a map on $L_2(SU_q(\ell+1)\verb1\1SU_q(\ell))$ by $\hat{\pi}$.
It is easy to check that $(\hat{\pi},\hat{u})$ is a covariant representation
for the system $(A_\ell,SU_q(\ell+1),\tau)$.

\bthm
 Spectral dimension of the odd dimensional quantum sphere $S_q^{2\ell+1}$ is $2\ell+1$.
\ethm
\prf
Let $D_{eq}$ be the operator on $L_2(S_q^{2\ell+1})$ given by:
\be\label{eq_sphere1}
D_{eq} e_{\bldr^{nk},\blds}=(n+k) e_{\bldr^{nk},\blds}.
\ee
It follows from theorem~6.4, \cite{cha-pal-2008a},  $(\cla(S_{q}^{2\ell+1}),L_2(S_q^{2\ell+1}),D_{eq})$  is an 
equivariant spectral triple.

Note that the eigenspace corresponding to the eigenvalue
$n$ is 
\[
\mbox{span\,}\{e_{\bldr^{n-k,k},\blds}: 0\leq k\leq n, \blds\in \Gamma^{n-k,k}_0\}
\].
Let $\delta_n$ denote the dimension of this space.

Observe that for a given $n$ and $k$, the set $\Gamma^{n-k,k}_0$ consists
of all  GT tableaux  of the form 
\[
\blds=\left(\begin{matrix}
  c_1=n& k  & k& \cdots& k & k & d_1=0 \cr
  c_2& k & k &\cdots & k & d_2&\cr
   \cdots    & &\cdots&&&\cr
   c_{\ell-1}& k& d_{\ell-1} &&&&\cr
   c_\ell & d_\ell &&&&&\cr
                     c_{\ell+1}=d_{\ell+1} &&&&&&
                   \end{matrix}\right)
\]
Therefore the number of GT tableaux $\blds\in \Gamma^{n-k,k}_0$
is  $O(n^{2\ell-1})$. Hence $\delta_n=O(n^{2\ell})$.
This implies that for $p>2\ell+1$, one has
\[
\mbox{Trace\,}|D_{eq}|^{-p}\leq 
\sum_n n^{-p}n^{2\ell}=\sum_n\frac{1}{n^{p-2\ell}}<\infty,
\]
i.e.\ $D_{eq}$ is $p$-summable.

Next, let us take an $\epsilon\in (0,\frac{1}{4\ell+2})$.
Then for large enough $n$, the number of Young tableux
$\lambda=\lambda_{n-k,k}\equiv(n,k,\ldots,k,0)$
with
\[
\left|k-\left(\frac{\ell+1}{2\ell+1}\right)n\right|<\epsilon n
\]
is of the order $n$.
For each such $\lambda_{n-k,k}$, the number of 
$\blds\in\Gamma^{n-k,k}_0$ with
\begin{multline*}
|d_2-\frac{1}{2\ell+1}n|<\epsilon n,\quad
 |d_3-\frac{2}{2\ell+1}n|<\epsilon n,\ldots,
  |d_{\ell+1}-\frac{\ell}{2\ell+1}n|<\epsilon n,\\
  |c_{\ell}-\frac{\ell+2}{2\ell+1}n|<\epsilon n,\quad
  |c_{\ell-1}-\frac{\ell+3}{2\ell+1}n|<\epsilon n,\ldots,
  |c_2-\frac{2\ell}{2\ell+1}n|<\epsilon n,
\end{multline*}
is $Cn^{2\ell-1}$, which implies that
$\delta_n\geq C n^{2\ell}$
(here $C$ denotes a generic constant, independent of $n$).
Hence for $p\leq 2\ell+1$, we have
\[
\mbox{Trace\,}|D_{eq}|^{-p}\geq \sum_n n^{-p}n^{2\ell}
  =\sum_n \frac{1}{n^{p-2\ell}}=\infty,
\]
i.e.\ $D_{eq}$ can not be $p$-summable.

Next, let $D$ be a self-adjoint operator with compact resolvent on $L_2(S_q^{2\ell+1})$
that is equivariant with respect to the covariant representation $(\hat{\pi},\hat{u})$.
By theorem~6.4, \cite{cha-pal-2008a}, one then has $|D|\leq a+bD_{eq}$ for some constants
$a$ and $b$. Therefore $D$ can not be $p$-summable for
$p\leq 2\ell+1$.

Thus the spectral dimension of the $C^*$-dynamical system under consideration is $2\ell+1$.
\qed

\section{$SU_q(2)$ action on Podle\'{s} sphere $S_{q0}^2$}
Quantum sphere was introduced by Podle\'s in \cite{pod-1987a}. 
The $C^*$-algebra $C(S_{q,0}^2)$ is the universal $C^*$-algebra generated 
by two elements $\xi$ and $\eta$ subject to the following relations:
\bean
\xi^* = \xi ,&&  \eta^*\eta = \xi - \xi^2,\\
\eta\xi = q^2 \xi\eta,&&\eta\eta^* = q^2 \xi - q ^ 4 \xi^2.
\eean
Here the deformation parameters $ q$ satisfies $ | q | < 1$. 
This space was studied in detail in \cite{pod-1987a}. Let us restate
the relevant facts from that paper in our present notation, so as to be able
to make use of the computations we have done in earlier sections. 

Let $u$ denote the fundamental unitary for $SU_q(2)$,
which is the irreducible unitary representation corresponding to the
Young tableaux $\one=(1,0)$.
Let $\bldz$ denote the function $t\mapsto t$ on the torus group $\bbt$.
Then $C(\bbt)$ is generated by the unitary $\bldz$.
Define $\phi:C(SU_q(2))\rightarrow C(\bbt)$ by
\be
\phi(u_{ij})=\begin{cases} \bldz & \mbox{if $i=j=1$},\cr
        \bldz^* & \mbox{if $i=j=2$},\cr
        0 & \mbox{if $i\neq j$}.
       \end{cases}
\ee
This is a quantum group homomorphism from $C(SU_q(2))$ onto $C(\bbt)$.
Then 
\[
C(SU_q(2)\verb1\1\bbt):=\{a\in C(SU_q(2): (\phi\otimes\id)\Delta(a)=I\otimes a\}
\]
 is the $C^*$-subalgebra of
$C(SU_q(2))$ generated by the elements $u_{11}u_{21}$ and $u_{12}u_{21}$.
Define $\psi:C(S_{q0}^{2})\rightarrow C(SU_q(2)\verb1\1\bbt)$
by 
\[
 \psi(\xi)=-q^{-1}u_{12}u_{21}, \quad \psi(\eta)=u_{11}u_{21}.
\]
This gives an isomorphism between $C(SU_q(2)\verb1\1\bbt)$ and
$C(S_{q0}^{2})$,
and the action $\tau$ of $SU_q(2)$ on $S_{q0}^2$ is the action induced from
the comultiplication map of $SU_q(2)$, i.e.\ the homomorphism
that makes the following diagram commute:
\[
\def\labelstyle{\scriptstyle}
  \xymatrix@C=23pt@R=50pt{
 C(S_{q0}^{2})\ar[d]_{\psi}\ar[r]_-{\tau} & C(S_{q0}^{2})\otimes
C(SU_q(2))\ar[d]_{\psi\otimes\id}&\\
C(SU_q(2)\backslash\bbt)\ar[r]_-{\Delta}  & C(SU_q(2)\backslash\bbt)\otimes C(SU_q(2))\\
}
\]
and  the invariant state for this action is
the restriction of the Haar state on $C(SU_q(2))$ to this $C^*$-subalgebra.

Let $A(S_{q0}^2)$ denote the involutive algebra generated by $\xi$ and $\eta$.
Then one has
\[
A(S_{q0}^2)=\mbox{span}\{e_{\bldr^{kk}\blds}: k\in\mathbb{N}, \blds\in\Gamma^{kk}_0\}
= \mbox{span}\{e_{\bldr^{kk}\bldr^{2k-m,m}}: k,m\in\mathbb{N}, 0\leq m\leq 2k\}.
\]

\bthm
The spectral dimension of $S_{q0}^2$ is 0.
\ethm
\prf
It follows from the above discussion that the $L_2$ space $L_2(S_{q0}^2)$ of the sphere is  the closed subspace of $L_2(SU_q(2))$ 
spanned by 
$\{e_{\bldr^{kk}\bldr^{2k-m,m}}: k,m\in\mathbb{N}, 0\leq m\leq 2k\}$
and the representation of $C(S_{q0}^2)$ on this
is  the restriction of the left multiplication representation $\pi$ of
$C(SU_q(2))$ to $C(S_{q0}^2)$. We will call this restriction $\hat{\pi}$.
From equation~(\ref{eq:left_mult}), it follows that the 
action of the elements $\xi$ and $\eta$ on
the basis elements are given by:
\bea
\lefteqn{\hat{\pi}(\xi)e_{\bldr^{kk}\bldr^{2k-m,m}} }\nonumber\\
&=& 
      - q^{k+m-1}e_{\bldr^{k-1,k-1}\bldr^{2k-m-1,m-1}} + (q^{2k}+q^{2m})e_{\bldr^{kk}\bldr^{2k-m,m}}
          - q^{k+m+1}e_{\bldr^{k+1,k+1}\bldr^{2k-m+1,m+1}} \nonumber\\
&&  \\
\lefteqn{\hat{\pi}(\eta)e_{\bldr^{kk}\bldr^{2k-m,m}} }\nonumber\\
&=& 
      - q^{k}e_{\bldr^{k-1,k-1}\bldr^{2k-m,m-2}} + q^m (1-q^{2k})e_{\bldr^{kk}\bldr^{2k-m+1,m-1}}
          - q^{k+2m+2}e_{\bldr^{k+1,k+1}\bldr^{2k-m+2,m}}.\nonumber\\
&& 
\eea
Clearly the subspace $L_2(S_{q0}^2)$ is an invariant subspace for the right regular representation
of $SU_q(2)$. Let us call this restriction $\hat{u}$. Then $(\hat{\pi},\hat{u})$ gives a covariant
representation for the system $(S_{q0}^2, SU_q(2),\tau)$ on $L_2(S_{q0}^2)$.
Restriction of $\hat{u}$ to $\mbox{span}\{e_{\bldr^{kk}\bldr^{2k-m,m}}: 0\leq m\leq 2k\}$
is equivalent to the irreducible $u^{(2k,0)}$. Therefore any equivariant Dirac operator $D$ will
be of the form
\[
e_{\bldr^{kk}\bldr^{2k-m,m}}\mapsto d(k)e_{\bldr^{kk}\bldr^{2k-m,m}},\quad k\in\bbn, 0\leq m \leq 2k.
\]
A necessary and sufficient condition for the boundedness of the commutators 
of $D$ with the $\hat{\pi}(\xi)$ and $\hat{\pi}(\eta)$ (and hence with all $\hat{\pi}(a)$'s) is the following:
\be\label{eq:bddness_2sphere}
|d(k)-d(k+1)| =O( q^{-k}).
\ee
Thus an operator $D$ given by
\[
e_{\bldr^{kk}\bldr^{2k-m,m}}\mapsto q^{-k}e_{\bldr^{kk}\bldr^{2k-m,m}},
               \quad k\in\bbn, 0\leq m \leq 2k
\]
makes $(L_2(S_{q0}^2), \hat{\pi},D)$ an equivariant spectral triple.
Dimension of the eigenspace corresponding to the eigenvalue
$q^{-k}$ is $2k+1$. Therefore it follows that for any $p>0$, this
spectral triple is $p$-summable. Thus the spectral dimension
of the Podle\'s sphere $S_{q0}^2$ is 0.\qed

\section{$A_u(Q)$ action on Cuntz algebras}
In this section we are going to compute the spectral dimension of the
Cuntz algebras. For a brief account on Cuntz algebras, see \cite{cun-1977a}.
Given a nonsingular $n \times n$ matrix $Q \in
GL(n,{\mathbb C})$, the universal quantum groups $A_u(Q)$ were
introduced by Van Daele and Wang in (\cite{van-wang-1996a}) as the
universal compact quantum group $(A_u(Q),u)$ generated by $u_{ij}, (1
\le i,j \le n)$ with defining relations
\bean
u^*u=I_n=uu^*,u^tQ\overline{u}Q^{-1}=I_n=Q\overline{u}Q^{-1}u^t,
\eean
where $u=(u_{ij}), {(\overline{u})}_{ij}=u_{ij}^*$ and
${(u^t)}_{ij}=u_{ji}$. Wang showed (\cite{wang-1999a}) that $A_u(Q)$, with $Q$ positive of
trace 1 acts ergodically on the Cuntz algebra ${\mathcal O}_n$. Recall 
that ${\mathcal O}_n$ is the universal $C^*$-algebra generated by
$n$-isometries $S_k(k=1,\cdots,n)$ such that $\sum S_kS_k^*=1$ and the
action is specified by $\tau(S_j)=\sum_{i=1}^n S_i \otimes u_{ij}$. The
dense $*$-subalgebra generated by the $S_i$'s gives the span of spectral
subspaces. Here is an explicit description of the spectral subspaces
following \cite{wang-1999a}. For a multi-index $\alpha=(i_1,\cdots,i_r), 1 \le
i_1,\cdots,i_r \le n$ of length $\ell(\alpha )=r$, let $S_\alpha=
S_{i_1}\cdots S_{i_r}$. Then the spectral subspaces are given by
${\mathcal H}_{rs}=\mbox{Span\,} \{ e_{\alpha, \beta}=S_\alpha
S_\beta^*: \ell(\alpha)=r,\ell(\beta)=s \}$. An equivariant ``Dirac'' operator
$D$ must be constant on these subspaces, hence they are given by
\[
D:e_{\alpha, \beta} \mapsto d(\ell(\alpha), \ell(\beta)) e_{\alpha,\beta},
\]
where $d: {\mathbb N} \times {\mathbb N} \rightarrow {\mathbb R}$ is a
real valued function. For a multi-index $\alpha=(i_1,\cdots,i_r)$ let
$S_i\alpha$ and $S_i^{-1}\alpha$ denote the multi-indices given by   
$S_i \alpha= (i,i_1,\cdots,i_r), S_i^{-1}\alpha =(i_2,\cdots,i_r)$, if
$i_1=i$ and $e_{S_i^{-1} \alpha,\beta}= 0$, if $i \ne i_1$. In this
notation we have $S_i e_{\alpha,\beta}=e_{S_i \alpha, \beta} $ and
$S_i^*e_{\alpha,\beta}=e_{S_i^{-1} \alpha \beta}$. Let $\rho_Q$ be the
unique invariant state and M be such that $\|[D,S_i]\|, \|[D,S_i^*]\| <
M$ for $1 \le i \le n$. Then
\bean
[D,S_i]e_{\alpha,\beta}& = &
(d(\ell(\alpha)+1,\ell(\beta))-d(\ell(\alpha),\ell(\beta)))e_{S_i \alpha,\beta}, \\
\rho_Q(e^*_{\alpha,\beta} e_{\alpha,\beta}) & = & \rho_Q(S_\beta
S^*_\alpha S_\alpha S^*_\beta)\\
& = & \rho_Q( S_\beta S^*_\beta) \\
& = & \rho_Q(e_{\beta, \beta}) \\
& = & \rho_Q(e^*_{S_i \alpha,\beta} e_{S_i \alpha,\beta)}\\
\eean
together imply 
\bea \label{cuntz1}
|d(\ell(\alpha)+1,\ell(\beta))-d(\ell(\alpha),\ell(\beta))| < M\quad \forall
\alpha,\beta .
\eea
Similarly using $S^*_i e_{\phi,\beta}=e_{\phi,\beta'}$, with 
$\beta'=(j_1,\cdots,j_s,i)$, where $\beta=(j_1,\cdots,j_s)$ and $$
[D,S^*_i]e_{\phi,\beta}= (d(0,\ell(\beta))-d(0,\ell(\beta)))e_{\phi,\beta} $$
we get 
\bea \label{cuntz2} |d(0,\ell(\beta))-d(0,\ell(\beta))| < M. \eea
Thus combining (\ref{cuntz1}) and (\ref{cuntz2}) we get
$|d(\ell(\alpha),\ell(\beta))| < M (\ell(\alpha) + \ell(\beta))$. Therefore $D$ has
eigenvalue $k$ with multiplicity same as the cardinality of $\{
e_{\alpha,\beta}: \ell(\alpha)+\ell(\beta)=k\}$, which is $n^k$. Clearly there
is no positive number $s$ such that $Tr |D|^{-s}= \sum_k n^k k^{-s} <
\infty$. Hence ${\mathcal Sdim}({\mathcal O}_n, A_u(Q), \tau)=\infty$.
Thus we have proved:
\begin{thm}
 Spectral dimension of Cuntz algebra is infinite.
\end{thm}

\section{Concluding remarks}
Few remarks are in order.

\begin{enumerate}
\item
In defining the invariant, we have not demanded nontriviality (of the
$K$-homology class) of the spectral triples that we take into consideration.
One reason behind that is the following. For many examples, if we just take
the GNS space of the invariant state to be the Hilbert space where the spectral
triples live and demand nontriviality, then there may not exist any such spectral triple.
This, for example, is the case for the entire family $SU_q(\ell+1)$, $\ell>1$
(see \cite{cha-pal-2005a}). On the other extreme, in the case of the
classical $SU(2)$ for instance, as theorem~(5.4) of (\cite{cha-pal-2003a})
shows, the spectral dimension of $SU(2)$ will come out to be 4.
On the other hand, if we want to allow Hilbert spaces other than
just one copy of the GNS space, that takes us back to one of the problems that was mentioned in
the introduction, namely, the choice of a `natural' Hilbert space.
One somewhat natural choice might be to look at the GNS space tensored with
$\mathbb{C}^n$ where one allows $n$ to be greater than one. But then another 
problem that we discussed in section~1 crops up, namely the collection
of such spectral triples becomes big and somewhat intractable.

\item
The fact that Podle\'{s} sphere 
has spectral dimension zero should not come as a surprise, 
because as shown in (\cite{nes-tus-2005a}), on many counts 
its behaviour differs from the classical case.
This particular example also illustrates that the invariant for
a homogeneous space of the $q$-deformation of a classical Lie group is not 
merely the geometric dimension of its classical counterpart.

\item
 Noncommutative Geometry offers a new way of looking at classical situations. For example,
 we have seen in section~3
 that for classical $SU(2)$ its spectral dimension is same as its dimension as a Lie group. 
 It would then be tempting to conjecture that 
 \begin{quotation}
  \noindent\textit{the spectral dimension of a homogeneous space of a (classical) compact Lie group is
    same as its dimension as a differentiable manifold.}
    \end{quotation}
  It must be pointed out however that at the moment, except for the $SU(2)$ example,
  nothing much is known that one can cite as a strong evidence. The case of
  quantum $SU(\ell+1)$ and   the quantum odd dimensional spheres do point towards
  the above statement, but one has to keep in mind that the behaviour of these spaces
  can be quite different when $q=1$.
  So this should perhaps be looked upon as more of  a pointer to future research.
  Also, if the above conjecture turns out to be false, then spectral dimension 
  produces a new invariant; the question then is: is it some classically known quantity 
  associated with the space?
\end{enumerate}

%

\noindent{\sc Partha Sarathi Chakraborty}
(\texttt{parthac@imsc.res.in})\\
         {\footnotesize  Institute of Mathematical Sciences, 
CIT Campus, Chennai--600\,113, INDIA}\\[1ex]
{\sc Arupkumar Pal} (\texttt{arup@isid.ac.in})\\
         {\footnotesize Indian Statistical
Institute, 7, SJSS Marg, New Delhi--110\,016, INDIA}

\end{document}